\input ./style/arxiv-general.cfg
\documentclass[aap,MSNbibl,nameyear,dvips]{arximspdf}
\makeatletter
   \@ifpackageloaded{graphicx}{}{\usepackage{graphicx}}
\makeatother
\usepackage{mathbh}

\doi{10.1214/14-AAP1044} 
\volume{25}
\issue{4}
\pubyear{2015}
\firstpage{2096}
\lastpage{2133}
\docsubty{FLA}

\makeatletter
\newcommand{\F}{\mathrm{F}}
\newcommand{\R}{\mathrm{R}}
\newcommand{\ext}{\mathrm{ext}}
\newcommand{\rrVert}{\Vert}
\newcommand{\llVert}{\Vert}
\newtheorem{teo}{Theorem}[section]

\newtheorem{lem}[teo]{Lemma}
\newtheorem{prop}[teo]{Proposition}
\newproclaim{defin}[teo]{Definition}
\newproclaim{rem}[teo]{Remark}
\newtheorem{theorem}[teo]{Theorem}
\newtheorem{lemma}[teo]{Lemma}
\newtheorem{proposition}[teo]{Proposition}

\makeatother

\begin{document}
\begin{frontmatter}

\title{Global solvability of a networked integrate-and-fire model of McKean--Vlasov type\thanksref{T1}}
\runtitle{An integrate-and-fire model of McKean--Vlasov type}

\begin{aug}
\author[A]{\fnms{Fran\c cois}~\snm{Delarue}\ead[label=e1]{Francois.Delarue@unice.fr}},
\author[B]{\fnms{James}~\snm{Inglis}\corref{}\ead[label=e2]{James.Inglis@inria.fr}},\\
\author[A]{\fnms{Sylvain}~\snm{Rubenthaler}\ead[label=e3]{Sylvain.Rubenthaler@unice.fr}}
\and
\author[B]{\fnms{Etienne}~\snm{Tanr\'e}\ead[label=e4]{Etienne.Tanre@inria.fr}}
\runauthor{Delarue, Inglis, Rubenthaler and Tanr\'e}
\affiliation{Universit\'e Nice Sophia Antipolis,
Inria Sophia Antipolis---M\'editerran\'ee,\\
Universit\'e Nice Sophia Antipolis
and
Inria Sophia Antipolis---M\'editerran\'ee}
\address[A]{F. Delarue\\
S. Rubenthaler\\
Laboratoire J. A. Dieudonn\'e\\
Universit\'e Nice Sophia Antipolis\\
Parc Valrose\\
06108 Nice Cedex 02\\
France\\
\printead{e1}\\
\phantom{E-mail: }\printead*{e3}}
\address[B]{J. Inglis\\
E. Tanr\'e\\
Equipe Tosca\\
Inria Sophia Antipolis---M\'editerran\'ee\\
2004 route des lucioles,
BP 93\\
06902 Sophia Antipolis Cedex\\
France\\
\printead{e2}\\
\phantom{E-mail: }\printead*{e4}}
\end{aug}
\thankstext{T1}{Supported by the Agence National de la Recherche
through the ANR Project MANDy ``Mathematical Analysis of Neuronal
Dynamics,'' ANR-09-BLAN-0008-01.}

\received{\smonth{4} \syear{2014}}

%
\begin{abstract}
We here investigate the well-posedness of a networked
integrate-and-fire model describing an infinite population of neurons
which interact with one another through their common statistical
distribution. The interaction is of the self-excitatory type as, at any
time, the potential of a neuron increases when some of the others fire:
precisely, the kick it receives is proportional to the instantaneous
proportion of firing neurons at the same time. From a mathematical
point of view, the coefficient of proportionality, denoted by $\alpha
$, is of great importance as the resulting system is known to blow-up
for large values of $\alpha$. In the current paper, we focus on the
complementary regime and prove that existence and uniqueness hold for
all time when $\alpha$ is small enough.
\end{abstract}

%
\end{frontmatter}

\setcounter{footnote}{1}

\section{Introduction}\label{sec1}

The stochastic integrate-and-fire model for the membrane potential $V$
across a neuron in the brain has received a huge amount of attention
since its introduction [see \citet{sacerdotegiraudo} for a
comprehensive review]. The central idea is to model $V$ by threshold
dynamics, in which the potential is described by a simple linear
(stochastic) differential equation up until it reaches a fixed
threshold value $V_\F$, when the neuron emits a ``spike''.
Experimentally, at this point an action potential is observed, whereby
the potential increases very rapidly to a peak (hyperpolarization
phase) before decreasing quickly to a reset value (depolarization phase).

Since spikes are stereotyped events, they are fully characterized by
the times at which they occur. The integrate-and-fire model is part of
a family of spiking neuron models which take advantage of this by
modeling only the spiking times and disregarding the nature of the
spike itself. Specifically, in the integrate-and-fire model we observe
jumps in the action potential as the voltage is immediately reset to a
value $V_\R$ whenever it reaches the threshold $V_\F$. Despite its
simplicity, versions of the integrate-and-fire model have been able to
predict the spiking times of a neuron with a reasonable degree of
accuracy [\citet{jolivetlewisgerstner,kistlergerstnervanhemmen}].

Many extensions of the basic integrate-and-fire model have been studied
in the neuroscientific literature, including ones in which attempts are
made to include noise and to describe the situation when many
integrate-and-fire neurons are placed in a network and interact with
each other. In \citet{lewisrinzel,ostojicbrunelhakim}, the
following equation describing how the potential $V_i$ of the $i$th
neuron in a network of $N$ behaves in time is proposed:
%
%
\begin{eqnarray}\label{Brunel}
\frac{d}{dt}V_i(t) &=& -\lambda
V_i(t) + \frac{\alpha}{N}\sum_{j}\sum
_{k} \delta_0\bigl(t-
\tau_{k}^j\bigr)
\nonumber\\[-8pt]\\[-8pt]
&&{} + \frac{\beta}{N}\sum_{j\neq i}V_j(t) + I^{\ext}_i(t) +\sigma\eta_i(t)\nonumber
\end{eqnarray}
for $V_i(t) < V_\F$ and where $V_i(t)$ is immediately reset to $V_\R$
when it reaches $V_\F$. Here, $I^{\ext}_i(t)$ represents the external
input current to the neuron, $\eta_i(t)$ is the noise (a white noise)
which is importantly supposed to be independent from neuron to neuron,
and the constants $\lambda, \beta, \alpha$ and $\sigma$ are chosen
according to experimental data. Moreover, the interaction term is
described in terms of $\tau_{k}^j$, which is the time of the $k$th
spike of neuron $j$, and the Dirac function $\delta_0$. Precisely, it
says that whenever one of the other neurons in the network spikes, the
potential across neuron $i$ receives a ``kick'' of size $\alpha/N$.
The Dirac mass interactions give rise to the same kind of instantaneous behavior
as the integrate-and-fire model. Although it is a simplification of
reality, it produces some interesting phenomena
from a biological perspective [see \citet{ostojicbrunelhakim}].

In the case of a large network, that is, when $N$ is large, many
authors approximate the interaction term
by an instantaneous rate $\nu(t)$, the so-called mean-firing rate
[see, e.g., \citet
{brunel,brunelhakim,ostojicbrunelhakim,renartbrunelwang}]. However, in
the neuroscience literature, little attention is paid to how this
convergence is achieved.
Mathematically, the mean-field limit as $N\to\infty$ must be taken,
but as a first step,
this requires a careful analysis of the asymptotic well-posedness.
This is precisely the purpose of the paper: to focus on the unique
solvability of the resulting nonlinear limit equation
(the analysis of the convergence being left to
further investigations). At first glance such a question may seem
classical, given the volume of results available that guarantee the
existence of a solution to distribution dependent SDEs. However, as
quickly became apparent in our analysis, in the excitatory case
($\alpha>0$) the problem is in fact a delicate one, for which, to our
knowledge, there are no existing results available. This difficulty is
due to the nature of the interactions, which introduce the strong
possibility of a solution that ``blows up'' in finite time. The
validity of the study of this question, and its nontrivial nature, is
further justified by the fact that several authors have recently been
interested in exactly the same problem from a PDE perspective
[\citet{perthame,carillogonzalesgualdanischonbek}]. Despite some
serious effort and very interesting related results on their part, we
understand that they were not able to prove the
existence and uniqueness of global solutions to the limit equation,
which is the main result of the present paper.

\subsection{Precisions}
We now make precise the nonlinear equation of interest. First,
since the mathematical difficulties lie within the jump interaction term,
we suppose that there is no external input current [$I^{\ext}_i(t)
\equiv0$], and that the interaction term is composed solely of the
jump or reset part ($\beta= 0$). Although this is a nontrivial
simplification from a neuroscience perspective, it still captures all
the mathematical complexity of the resulting mean-field equation.

Without loss of generality, we also take the firing threshold $V_\F =1$
and the reset value $V_\R=0$ for notational simplicity. The
nonlinear stochastic mean-field equation under study here is then
%
%
\begin{equation}
X_{t}=X_0 + \int_{0}^{t}b(X_{s})\,ds+
\alpha\mathbb{E} (M_{t})+\sigma W_{t}-M_{t},
\qquad t\geq0,\label{Intrononlineareq}
\end{equation}
where $X_0<1$ almost surely, $\alpha\in\mathbb{R}$, $\sigma>0$,
$(W_{t})_{t\geq0}$ is a standard Brownian motion in $\mathbb{R}$ and
$b\dvtx \mathbb{R}\to\mathbb{R}$ is Lipschitz continuous. In comparison with
(\ref{Brunel}), $b$ must be thought of as $b(x)=-\lambda x$.
Equation (\ref{Intrononlineareq}) is then intended to describe the
potential of one \textit{typical} neuron in the infinite network, its
jumps (or resets) being given by
\[
M_{t}=\sum_{k\geq1}\mathbh{1}_{[0,t]}(
\tau_{k}),
\]
where $(\tau_{k})_{k\geq1}$ stands for the sequence of hitting times
of $1$ by the process $(X_t)_{t\geq0}$.
That is, $(M_{t})_{t\geq0}$
counts the number of times $X_{t}$ hits the threshold before time $t$,
so that $\mathbb{E}(M_t)$ denotes the
\textit{theoretical}
expected number of times the threshold is reached before $t$.
Such a theoretical expectation corresponds to what we would envisage as
the limit of the integral form of the interaction term
\[
\frac{1}{N} \int_{0}^t \sum
_{j}\sum_{k} \delta\bigl(s-
\tau_{k}^j\bigr) \,ds = \frac{1}{N} \sum
_{j} \sum_{k} {\mathbh 1}_{\{\tau_{k}^j \leq t\}}
\]
in (\ref{Brunel}) when $N\to\infty$, assuming that neurons become
asymptotically independent [as is observed in more classical particle
systems---see \citet{sznitman}].

\subsection{PDE viewpoint and ``blow-up'' phenomenon} As mentioned
above, equation (\ref{Intrononlineareq}) has been rigorously studied
from the PDE viewpoint before. When $\sigma\equiv1$, the
Fokker--Planck equation for the density $p(t, y)\,dy = \mathbb{P}(X_t
\in dy)$ is given by
\[
\partial_t p(t, y) + \partial_y \bigl[ \bigl(b(y) +
\alpha e'(t) \bigr)p(t, y) \bigr] - \tfrac{1}{2}
\partial^2_{yy}p(t, y) = \delta_0(y)e'(t),
\qquad y< 1,
\]
where $e(t) = \mathbb{E}(M_t)$, subject to $p(t, 1) = 0$, $p(t,
-\infty) = 0$, $p(0, y)\,dy = \mathbb{P}(X_0\in dy)$.
Moreover, the condition that $p(t,y)$ must remain a probability density
translates into the fact that
\[
e'(t) = \frac{d}{dt}\mathbb{E}(M_t) = -
\frac{1}{2}\partial_yp(t,1)\qquad\forall t>0,
\]
which describes the nonlinearity of the problem.
In the case when $b(x) = -\lambda x$, this nonlinear Fokker--Planck
equation is exactly the one studied in \citet{perthame} and
\citet{carillogonzalesgualdanischonbek}. Therein, the authors
conclude that for some choices of parameters, no global-in-time
solutions exist. The term ``blow-up'' is then used to describe the
situation where the solution (defined in a weak sense) ceases to exist
after some finite time.
With our formulation, since $e'(t)$ corresponds to the mean firing rate
of the infinite network,
it is very natural to define a ``blow-up'' time as a time when $e'(t)$
becomes infinite. Intuitively, this can be understood as a point in
time at which a large proportion of the neurons in the network all
spike at exactly the same time, that is, the network \textit{synchronizes}.

In \citet{perthame} and \citet
{carillogonzalesgualdanischonbek}, it is shown that, in the cases
$\alpha=0$ and $\alpha< 0$ (the latter one being referred to as
``self-inhibitory'' in
neuroscience), the nonlinear Fokker--Planck equation has a unique
solution that does not blow-up in finite time.
However, in the so-called ``self-excitatory'' framework, that is, for
$\alpha>0$,
existence of a solution for all time is left open. Instead, a negative
result is established [\citet{perthame}, Theorem~2.2],
stating that, for any $\alpha>0$, it is possible to find an initial
probability distribution $\mathbb{P}(X_0 \in dy)$ such that any
solution must blow-up in finite time, that is, such that $e'(t) =
\infty$ for some $t>0$.

\subsection{Present contribution}
In this paper, we thus investigate the case $\alpha\in(0,1)$. Our
main contribution is to show that, given a starting point $X_0 = x_0$,
we can find an explicit $\alpha$ small enough so that there does
indeed exist a unique global-in-time solution to (\ref
{Intrononlineareq}) (and hence to the associated Fokker--Planck
equation) which does not blow-up (see Theorem~\ref{solutionuptoT}). In
view of the above discussions, our result complements and goes further
than those found in \citet{perthame} and \citet
{carillogonzalesgualdanischonbek}, and the surprising difficulty of the
problem is reflected in the rather involved nature of our proofs.

As already said, equation (\ref{Intrononlineareq}) can be thought of
as of McKean--Vlasov-type, since the process $(X_t)_{t \geq0}$
depends on the distribution of the solution itself. However, it is
highly nonstandard, since it actually depends on the distribution of
the \textit{first hitting times} of the threshold by the solution.
This renders the traditional approaches to McKean--Vlasov equations and
propagation of chaos, such as those presented in \citet
{sznitman}, inapplicable, because we have no {a priori}
smoothness on the law of the first hitting times. Thus, our results are
also new in this context.

The general structure of the proof is at the intersection between
probability and PDEs, the deep core
of the strategy being probabilistic. The main ideas are inspired from
the methods used to investigate the well-posedness of Markovian
stochastic differential equations involving some nontrivial
nonlinearity. Precisely, the first point is to tackle unique
solvability in small time: when the parameter $\alpha$ is (strictly)
less than 1 and the
density of the initial condition decays linearly at the threshold, it
is proved that the system induces a natural contraction in a
well-chosen space provided the time duration is small enough. In this
framework, the specific notion of a solution plays a crucial role as it
defines the right space for the contraction.
Below, solutions are sought in such a way that the mapping $e\dvtx  t
\mapsto{\mathbb E}(M_t)$ is continuously differentiable: this
is a crucial point as it permits to handle the process $(X_t)_{t \geq
0}$ as a drifted Brownian motion.
The second stage is then to extend existence and uniqueness from short
to long times. The point is to prove that some key quantity is
preserved as time goes by. Here, we prove that the system cannot
accumulate too much mass in the vicinity of 1. Equivalently, this
amounts to showing that the Lipschitz constant of the mapping $e\dvtx  t
\mapsto{\mathbb E}(M_t)$ cannot blow-up in a finite time. This is
where the condition $\alpha$ small enough comes in: when $\alpha$ is
small enough, we manage to give some estimates for the density of
$X_{t}$ in the neighborhood of $1$, the critical value of $\alpha$
explicitly depending upon the available bound of the density. Generally
speaking, we make use of standard Gaussian estimates of Aronson type
for the density. Unfortunately,
the estimates we use are rather poor as they mostly neglect the right
behavior of the density of $X_{t}$ at the boundary, thus yielding a
nonoptimal value. Anyhow, they serve as a starting point for proving a
refined estimate of the gradient of the density at the boundary: this
is the required ingredient for proving that, at any time $t$, the mass
of $X_{t}$ decays linearly in the neighborhood of 1, uniformly in $t$ in
compact sets, and thus to apply iteratively the existence and
uniqueness argument in small time. In this way,
we prove by induction that existence and uniqueness hold on any finite
interval and thus on the whole of $[0,\infty)$.

It is worth mentioning that the main lines for proving the
{a priori} estimate on the Lipschitz constant of $e\dvtx  t \mapsto
{\mathbb E}(M_t)$
are probabilistic, thus justifying the
use~of a stochastic approach to handle the model. Indeed, the key step
in the control of~the Lipschitz constant of
$e$ is an \textit{intermediate} estimate of H\"older type, the proof
of which is inspired from
the probabilistic arguments used in \citet{krylovsafonov} for
establishing the H\"older regularity
of solutions to nonsmooth PDEs.
\subsection{Prospects}
Our result is for a general Lipschitz function $b$, but there are two
important specific cases that we keep in mind: the Brownian case when
$b\equiv0$ and the Ornstein--Uhlenbeck case when $b(x) = -\lambda x$,
$\lambda\geq0$. The Ornstein--Uhlenbeck case is most relevant to
neuroscience, but surprising difficulties remain in the purely Brownian case.
In both of these cases, we are able to give an explicit $\alpha_0$
depending on the deterministic starting point $x_0$ such that (\ref
{Intrononlineareq}) has a global solution for all $\alpha<\alpha_0$.
However, our explicit values do not appear to be optimal: simulations
suggest that for a given $x_0$ there exist solutions that do not blow-up for $\alpha$ bigger than our explicit $\alpha_0$, while there
exist solutions that blow-up that do not satisfy the conditions of
\citet{perthame}. Thus, an interesting question is to determine
for a given initial condition the critical value $\alpha_c$ such that
for $\alpha<\alpha_c$ (\ref{Intrononlineareq}) does not exhibit blow-up.

Another point is to relax the notion of solution in order to allow the
mapping $e\dvtx
t \mapsto{\mathbb E}(M_{t})$ to be nondifferentiable (and thus to blow-up). From the modeling point of view, this would permit the description
of \textit{synchronization} in the network.
Actually, based on our understanding of the problem and numerical simulations,
our guess is that, in full generality, the mapping $e$ may be
decomposed into a sequence
of continuously differentiable pieces separated by isolated
discontinuities. In that perspective, we feel that our work could serve
as a basis for investigating the unique
solvability of solutions that blow-up. In order to design a proper
uniqueness theory, it seems indeed quite mandatory to
understand how general solutions behave in the continuously
differentiable regime (which is the precise
purpose of the present paper), and then how discontinuities can emerge
(which is left to further works).

The layout of the paper is as follows. We present the main results in
Section~\ref{secMainresults}.
Solutions are defined in Section~\ref{Solutionasafixedpoint}
while
Section~\ref{sectionExistenceanduniquenessinsmalltime} is devoted to
proving the existence and uniqueness
in small time. The proof of Theorem~\ref{teogradientbd} is given in
Section~\ref{Long-timeestimates}.

\section{Main results}\label{secMainresults}

\subsection{Set-up} As stated in the \hyperref[sec1]{Introduction}, we are interested
in solutions to the nonlinear McKean--Vlasov-type SDE
%
%
\begin{equation}
X_{t}=X_0 + \int_{0}^{t}b(X_{s})\,ds+
\alpha\mathbb{E}(M_{t})+W_{t}-M_{t},\qquad t
\geq0,\label{simplifiedeq}
\end{equation}
where $X_0<1$ almost surely, $\alpha\in(0,1)$ and $(W_{t})_{t \geq
0}$ is a standard Brownian motion
with respect to a filtration $({\mathcal F}_{t})_{t \geq0}$ satisfying
the usual conditions. The jumps, or resets, of the system are described by $(\tau_{0}=0)$
%
%
\begin{equation}\label{M}
\qquad M_{t}=\sum_{k\geq1}
\mathbh{1}_{[0,t]}(\tau_{k})\qquad\mbox{with } \tau_{k}
= \inf\{t>\tau_{k-1}\dvtx  X_{t-} \geq1\}, k\geq1.
\end{equation}
We assume that $b\dvtx (-\infty, 1]\to\mathbb{R}$ is Lipschitz continuous
such that
\[
\bigl|b(x)\bigr|\leq\Lambda\bigl(|x|+1\bigr), \qquad\bigl|b(x) - b(y)\bigr| \leq K|x-y| \qquad
\forall x,
y \in(-\infty, 1].
\]

%
\begin{rem}
\label{remsigma}
By the time change $u = t/\sigma^2$, we could handle more general
cases when the intensity of the noise
in (\ref{simplifiedeq}) is $\sigma>0$ instead of $1$.
\end{rem}

As discussed in the \hyperref[sec1]{Introduction}, the key point is to look for a
solution for which $t\mapsto\mathbb{E}(M_t)$ is continuously
differentiable, which would correspond to a solution that does not
exhibit a finite time blow-up. This leads to the following definition
of a solution to (\ref{simplifiedeq}), where as usual
$\mathcal{C}^1[0, T]$ denotes the space of continuously differentiable
functions on $[0, T]$.

\begin{defin}[{[Solution to (\ref{simplifiedeq})]}]\label{definitionsolution}
The process $(X_t, M_t)_{0\leq t \leq T}$ will be said to be a solution
to (\ref{simplifiedeq}) up until time $T$ if
$(M_t)_{0\leq t \leq T}$ satisfies (\ref{M}), the map $([0,T] \ni t
\mapsto\mathbb{E}(M_t))\in\mathcal{C}^1[0, T]$ and $(X_t)_{0\leq t
\leq T}$ is a strong solution of (\ref{simplifiedeq}) up until time $T$.
\end{defin}

\subsection{Statements}Our main result is given by the following two
theorems. The first guarantees that, when $\alpha$ is small enough, if
there exists a solution to (\ref{simplifiedeq}) on some finite time
interval, then the solution does not blow-up on this interval.

%
\begin{theorem}
\label{teogradientbd}
For a given $\varepsilon\in(0,1)$, there exists a positive constant
$\alpha_{0} \in(0,1]$, depending only upon $\varepsilon$, $K$ and $
\Lambda$, such that,
for any $\alpha\in(0, \alpha_{0})$ and any positive time
$T>0$, there exists a constant ${\mathcal M}_{T}$, only depending on
$T$, $\varepsilon$, $K$ and $\Lambda$, such that, for any initial
condition $X_{0}=x_{0} \leq1 - \varepsilon$, any solution to (\ref
{simplifiedeq}) according to Definition~\ref{definitionsolution}
satisfies $(d/dt)\mathbb{E}(M_{t}) \leq{\mathcal M}_{T}$, for all $t
\in[0,T]$.
\end{theorem}

The second theorem is the main global existence and uniqueness result.

%
\begin{teo}
\label{solutionuptoT}
For any initial condition $X_{0}=x_{0} <1$ and $\alpha\in(0,\alpha
_0)$, where $\alpha_0 = \alpha_0(x_0)$ is as in Theorem~\ref
{teogradientbd} (taking $\varepsilon= 1 - x_0$),
there exists a unique solution to the nonlinear equation (\ref
{simplifiedeq}) on any $[0, T]$, $T>0$, according to Definition~\ref
{definitionsolution}.
\end{teo}

The size of the parameter $\alpha_0$ in Theorem~\ref{teogradientbd}
is found explicitly in terms of $\varepsilon, K$ and $\Lambda$
(Proposition~\ref{propholderbd1}), but more precisely it derives from
the fact that in the course of our proof we
must first show that, {a priori}, any solution on $[0, T]$ to
the nonlinear equation (\ref{simplifiedeq}) with $X_0 = x_0 \leq
1-\varepsilon$ satisfies\footnote{In the whole paper, we use the very
convenient notation $\frac{1}{dx}\mathbb{P}(X \in dx)$ to denote the
density at point $x$ of the random variable $X$ (whenever it exists). }
%
%
\begin{equation}
\label{keyfact} \frac{1}{dx}\mathbb{P}(X_t \in dx) <
\frac{1}{\alpha}, \qquad t\in[0, T],
\end{equation}
in a neighborhood of the threshold $1$ (see Lemma~\ref{lemholderbd5}).
It is this restriction that determines the $\alpha_0$ in Theorem~\ref
{teogradientbd}, so that it depends only on the best {a priori}
estimates available for the density on the left-hand side of (\ref
{keyfact}). The stated explicit choice for $\alpha_0$ in Proposition
\ref{propholderbd1} merely ensures that (\ref{keyfact}) holds for all
$\alpha<\alpha_0$ for any potential solution.

\subsection{Illustration: The Brownian case}
To further highlight the criticality of the system, we here
illustrate the blow-up phenomenon in the Brownian case.
Consider equation (\ref{simplifiedeq}) with $b\equiv0$, set $e(t) =
\mathbb{E}(M_t)$ and fix $X_0 = x_0<1$. Then the conditions of Theorem
\ref{solutionuptoT} are trivially satisfied, and so we know that there
exists a global-in-time solution for all $\alpha\in(0, \alpha_0(x_0))$.

\begin{figure}[t]

\includegraphics{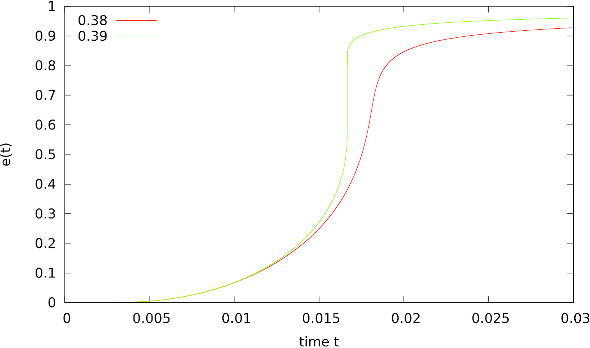}

\caption{Plot of $t\mapsto e(t)$ for $x_0=0.8$, $b(x)\equiv0$,
$\alpha=0.38$ \textup{(red)} and $\alpha=0.39$ \textup{(green)}.}
\label{figx008explosion}
\end{figure}

One may then ask if we ever observe a blow-up phenomenon in this case.
The affirmative answer can be seen by adapting the strategy in
\citet{perthame} [note that the result in \citet{perthame}
is written for
an Ornstein--Uhlenbeck type drift
but a similar argument applies when there is no drift].
For instance, choosing $x_{0}=0.8$,
computations show that global in time solvability must fail for $\alpha
\geq0.539$.
Moreover, tracking all the constants in the proof of Theorem
\ref{teogradientbd} below, we can find that $\alpha_0(0.8) \approx
0.104$, which
suggests that the system's behavior changes radically between these two values.
Such a radical change can be observed numerically
by investigating the graphs of $e(t) = \mathbb{E}(M_{t})$ for
different values of $\alpha$
in order to detect the emergence of some discontinuity.
Using a particle method to solve the nonlinear equation with $b\equiv
0$, we numerically observe
in Figure~\ref{figx008explosion}
that
the graph of $e$ is regular for $\alpha=0.38$ but has a jump
for $\alpha=0.39$. From the observations we have for other values of
$\alpha$, it seems that global solvability fails for $\alpha\geq
0.39$ and holds for $\alpha\leq0.38$.
%

As a summary, we present in Figure~\ref{fig6} the various regions of
the $\alpha$-parameter space $(0, 1)$ for $x_0 = 0.8$.
The region $\mathbf{D}$ stands for the set of $\alpha$'s for which
global solvability fails.
By the numerical experiments, it seems that global solvability also
fails in region $\mathbf{C}$, while by the same experiments it seems
that global solutions do exist for $\alpha$ in region $\mathbf{B}$.
In this article, we prove that
global solutions exist for~$\alpha\in\mathbf{A}$.\vspace*{-15pt}

\setlength{\unitlength}{1cm}
\begin{center}
%
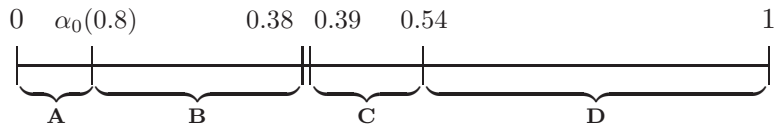
\begin{figure}[b]
\begin{picture}(10,2.3)
\label{summaryx0=08}
\put(0, 1.5){\line(1, 0){10}}
\put(0,1.25){\line(0, 1){.5}}
\put(-0.1, 2){$0$}
\put(10,1.25){\line(0, 1){.5}}
\put(9.9, 2){$1$}
\put(5.4,1.25){\line(0, 1){.5}}
\put(5.1,2){{\small$0.54$}}
\put(5.45, 1.25){$\underbrace{\hspace{4.5cm}}_\mathbf{D}$}
\put(3.9,1.25){\line(0, 1){.5}}
\put(3.95,2){{\small$0.39$}}
\put(3.95, 1.25){$\underbrace{\hspace{1.4cm}}_\mathbf{C}$}
\put(3.8,1.25){\line(0, 1){.5}}
\put(3.05,2){{\small$0.38$}}
\put(1.05, 1.25){$\underbrace{\hspace{2.7cm}}_\mathbf{B}$}
\put(1.0,1.25){\line(0, 1){.5}}
\put(0.5,2){{\small$\alpha_0(0.8)$}}
\put(0.05, 1.25){$\underbrace{\hspace{0.95cm}}_\mathbf{A}$}
\end{picture}\vspace*{-20pt}
\caption{Critical regions of $\alpha\in(0,1)$, for $x_0 = 0.8$ and
$b(x) \equiv0$.}
\label{fig6}
\end{figure}
\end{center}

\section{Solution as a fixed point}\label{Solutionasafixedpoint}

In this section, we identify a solution to the nonlinear equation (\ref
{simplifiedeq}) as a fixed point of an appropriate map on an
appropriate space. This will reduce the problem of finding a solution
to identifying a fixed point of this map.

Let $T>0$. For a general function $e\in\mathcal{C}^1[0, T]$, consider
the linear SDE
%
%
\begin{equation}
X_{t}^{e}=X_{0}+\int_{0}^{t}b
\bigl(X_{s}^{e}\bigr)\,ds+\alpha e(t)+W_{t}-M_{t}^{e},
\qquad t\in[0, T], X_{0}<1 \mbox{ a.s.},\hspace*{-20pt}\label{gammadefinition2}
\end{equation}
where $(W_{t})_{t \geq0}$ is a standard Brownian motion, $\alpha\in
(0, 1)$,
%
%
\begin{equation}
\label{defMt} M_{t}^{e}=\sum
_{k\geq1}\mathbh{1}_{[0,t]}\bigl(\tau_{k}^{e}
\bigr)
\end{equation}
and $\tau_{k}^{e}=\inf\{t>\tau_{k-1}^{e}\dvtx X_{t-}^{e}\geq1\}$ for
$k\geq1$,
$\tau_{0}^{e}=0$. The drift function $b$ is assumed to be Lipschitz as
above. Note that the solution to this SDE is well defined (by solving
(\ref{gammadefinition2}) iteratively from any
$\tau_{k}^e$ to the next $\tau_{k+1}^e$ and by noticing that the
jumping times $(\tau^e_{k})_{k \geq0}$ cannot accumulate in finite
time as the variations of $(X_{t}^e)_{t \geq0}$ on any $[\tau
_{k}^e,\tau_{k+1}^e)$, $k \geq0$, are controlled in probability). We
then define the map $\Gamma$ by setting
%
%
\begin{equation}
\Gamma(e) (t):=\mathbb{E}\bigl(M_{t}^{e}
\bigr)\label{gammadefinition}.
\end{equation}
We note that any fixed point of $\Gamma$ that is continuously
differentiable provides a solution to the nonlinear equation according
to Definition~\ref{definitionsolution} and vice versa. Thus, it is
natural to look for a fixed point of $\Gamma$ in a subspace of
$\mathcal{C}^1[0, T]$ where we are careful to uniformly control the
size of the derivative. Moreover, since it is clear from the definition
that $\Gamma(e)(0) = 0$ and $t\mapsto\Gamma(e)(t)$ is nondecreasing
for any $e \in\mathcal{C}^1[0, T]$, we in fact restrict the domain of
$\Gamma$ to the closed subspace $\mathcal{L}(T,A)$ of $\mathcal
{C}^1[0, T]$ defined by
\[
\mathcal{L}(T, A):= \Bigl\{ e\in\mathcal{C}^{1}[0,T]\dvtx  e(0)=0, e(s)
\leq e(t)\ \forall s\leq t, \sup_{0\leq t\leq T}e'(t)\leq A
\Bigr\}
\]
for some $A\geq0$. The map $\Gamma$ is thus defined as a map from
${\mathcal L}(T,A)$ into the set of nondecreasing functions on $[0,T]$.
It in fact depends on $A$ as its domain of definition depends on $A$;
for this reason, it should be denoted by $\Gamma^A$. Anyhow, since the
family $(\Gamma^A)_{A \geq0}$ is consistent in the sense that, for
any $A' \leq A$, the restriction of $\Gamma^A$ to ${\mathcal L}(T,A')$
coincides with $\Gamma^{A'}$, we can
use the simpler notation $\Gamma$.

The following {a priori} stability result provides further
information about where to look for fixed points, the proof of which we
leave until the end of the section.

%
\begin{prop} \label{stability} Given $T>0$, $a >0$ and $e \in
{\mathcal L}(T,A)$ it holds that
\begin{eqnarray*}
\label{stabilityequation}
&& \bigl( \bigl( \forall t \in[0,T], e(t)
\leq g_{a}(t)
\bigr)\mbox{ and }\bigl( \mathbb{E} \bigl[(X_{0})_{+} \bigr]
\leq a \bigr) \bigr)
\\
&&\qquad\Rightarrow \bigl( \forall t \in[0,T], \Gamma(e) (t) \leq g_{a}(t)
\bigr),
\end{eqnarray*}
where $(x)_+$ denotes the positive part of $x\in\mathbb{R}$, with
%
%
\begin{equation}
\label{gt} g_{a}(t):= \frac{a + (4 + \Lambda T^{1/2}) t^{1/2}}{1 -
\alpha} \exp\biggl(
\frac{2 \Lambda t}{1-\alpha} \biggr).
\end{equation}
\end{prop}

Letting $g(t):= g_{1}(t)$, $t \geq0$, since $X_0<1$ a.s., it thus
makes sense to look for fixed points of $\Gamma$ in the space
%
%
\begin{equation}
\label{H} \mathcal{H}(T, A):= \bigl\{ e\in\mathcal{L}(T, A)\dvtx  e(t)
\leq
g(t) \bigr\}.
\end{equation}
We equip $\mathcal{H}(T, A)$ with the norm $\|e\|_{\mathcal{H}(T,
A)}=\|e\|_{\infty,T}+\|e'\|_{\infty,T}$ inherited from $\mathcal
{C}^1[0, T]$. Here and throughout the paper, $\|\cdot\|_{\infty, T}$
denotes the supremum norm on $[0, T]$. $\mathcal{H}(T, A)$ is then a
complete metric space, since it is a closed subspace of~$\mathcal{C}^{1}[0,T]$.

For $e\in\mathcal{H}(T, A)$ Proposition~\ref{stability} implies that
$\Gamma(e)$ is finite and cannot grow faster that $g$, though it
remains to show that $\Gamma(e)$ is differentiable and that its
derivative is bounded by $A$
in order to check that $\Gamma$ indeed maps ${\mathcal H}(T,A)$ into
itself, for a suitable value of $A$ and $T$. The stability of
${\mathcal H}(T,A)$ by $\Gamma$ is discussed in Section~\ref{subse53}.

\subsection{Proof of Proposition \texorpdfstring{\protect\ref{stability}}{3.1}}
%
Fix $T>0$. We first note that we may write
%
%
\begin{eqnarray}
\label{integerrepresentation} M_{t}^{e}&=& \sup_{s \leq t}
\bigl\lfloor\bigl(Z_s^e \bigr)_{+} \bigr
\rfloor,
\nonumber
\\[-10pt]
\\[-10pt]
\nonumber
Z_{t}^e &=& X_{t}^{e}+M_{t}^{e}
= X_{0} + \int_{0}^t b
\bigl(X_{s}^e\bigr) \,ds + \alpha e(t) + W_{t},
\qquad t \in[0,T],
\end{eqnarray}
where $\lfloor x \rfloor$ denotes the floor part of $x\in\mathbb{R}$.
Indeed, one can see that for $t\in[\tau^e_{k},\tau^e_{k+1})$, $k\geq0$,
\begin{eqnarray*}
\sup_{s\leq t} \bigl\lfloor\bigl(Z_{s}^{e}
\bigr)_{+} \bigr\rfloor&=& \max\Bigl( \max_{0\leq j\leq k-1}
\Bigl( \sup_{s\in[\tau^e_{j},\tau^e_{j+1})} \bigl\lfloor\bigl
(X_{s}^{e}+j
\bigr)_{+} \bigr\rfloor\Bigr),\sup_{s\in[\tau^e_{k},t)} \bigl
\lfloor\bigl(X_{s}^{e}+k \bigr)_{+} \bigr\rfloor
\Bigr)
\\
&=&\max\Bigl( \max_{0\leq j\leq k-1} (j+1),k \Bigr) =M_{t}^{e},
\end{eqnarray*}
using the fact that $X_t^e<1$ for all $t\geq0$.

Then, given $t \in[0,T]$ such that $Z_{t}^e \geq0$, let $\rho^e:=
\sup\{ s \in[0,t]\dvtx  Z_{s}^e < 0\}$ ($\sup\varnothing= 0$). Pay
attention that $\rho^e$ is not a stopping time and that it depends on~$t$. Then, for $s \in[\rho^e,t]$,
%
%
\begin{equation}
\label{eq1291} \bigl| b\bigl(X_{s}^e\bigr) \bigr|\leq
\Lambda\bigl( 1 + \bigl| X_{s}^e \bigr|\bigr) \leq\Lambda
\bigl( 1 + \bigl| Z_{s}^e \bigr|+ M_{s}^e
\bigr) = \Lambda\bigl( 1 + \bigl(Z_{s}^e
\bigr)_{+} + M_{s}^e \bigr).
\end{equation}
By (\ref{integerrepresentation}), we know that $M_{s}^e \leq\sup_{0
\leq r \leq s} (Z_{r}^e)_{+}$. Therefore,
\[
\bigl| b\bigl(X_{s}^e\bigr) \bigr|\leq\Lambda\Bigl( 1+ 2
\sup_{0 \leq r \leq s} \bigl(Z_{r}^e
\bigr)_{+} \Bigr).
\]
By (\ref{integerrepresentation}), we obtain
%
%
\begin{equation}
\label{eq6101} Z_{t}^e \leq Z_{\rho^e}^e
+ \Lambda\int_{\rho^e}^{t} \Bigl( 1 + 2 \sup
_{0 \leq r \leq s} \bigl(Z_{r}^e\bigr)_{+}
\Bigr) \,ds + \alpha e(t) + W_{t} - W_{\rho^e}.
\end{equation}
If $\rho^e>0$, then $Z_{\rho^e}^e =0$ as, obviously, $(Z_{s}^e)_{0
\leq s \leq T}$ is a continuous process. If \mbox{$\rho^e=0$,}
then $X_{0} = Z_{\rho^e}^e \geq0$ since $Z_{\rho^e}^e$ is nonnegative.
Therefore,
%
%
\begin{equation}
\label{eq942} \qquad\bigl(Z_{t}^e\bigr)_{+} \leq
(X_{0})_{+} + \Lambda\int_{0}^{t}
\Bigl( 1 + 2 \sup_{0 \leq r \leq s} \bigl(Z_{r}^e
\bigr)_{+} \Bigr) \,ds + \alpha e(t) + 2 \sup_{0 \leq s \leq t}
| W_{s} |.
\end{equation}
Obviously, the above inequality still holds if $Z_{t}^e \leq0$.
We then notice that the process $(\sup_{0 \leq r \leq t}
(Z_{r}^e)_{+})_{0 \leq t \leq T}$
has finite values as $(Z^e_{t})_{0 \leq t \leq T}$ is continuous.
Therefore, taking the supremum in the left-hand side, applying
Gronwall's lemma and taking the expectation,
we deduce that $\mathbb{E} [ \sup_{0 \leq t \leq T}
(Z_{t}^e)_{+}]$ is finite. Taking directly the expectation in (\ref
{eq942}), we see that
%
%
\begin{eqnarray}
\label{eq1991} && \mathbb{E} \Bigl[ \sup_{0 \leq s \leq t}
\bigl(Z_{s}^e\bigr)_{+} \Bigr]
\nonumber\\[-8pt]\\[-8pt]\nonumber
&&\qquad \leq
\mathbb{E} \bigl[(X_{0})_{+}\bigr] + \Lambda\int
_{0}^{t} \Bigl( 1 + 2 \mathbb{E} \Bigl[ \sup
_{0 \leq r \leq s} \bigl(Z_{r}^e\bigr)_{+}
\Bigr] \Bigr) \,ds + \alpha e(t) + 4 t^{1/2},
\end{eqnarray}
for all $t\in[0, T]$. In particular, if ${\mathbb E}[(X_{0})_{+}] \leq
a$, $e(t) \leq g_{a}(t)$ for all $ t \in[0,T]$ [where $g_a$ is given
by (\ref{gt})],
and $R^{e}$ is the deterministic hitting time
\[
R^{e}:= \inf\Bigl\{ t \in[0,T]\dvtx  \mathbb{E} \Bigl[ \sup
_{0 \leq s \leq t} \bigl(Z_{s}^e\bigr)_{+}
\Bigr] > g_{a}(t) \Bigr\} \qquad(\inf\varnothing= + \infty),
\]
then, for any $t \in(0,R^e \wedge T]$,
\begin{eqnarray*}
&& \mathbb{E} \Bigl[ \sup_{0 \leq s \leq t} \bigl(Z_{s}^e
\bigr)_{+} \Bigr]
\\
&&\qquad  \leq a + \Lambda\int_{0}^{t}
\bigl( 1 + 2 g_{a}(s) \bigr) \,ds + \alpha g_{a}(t) + 4
t^{1/2}
\\
&&\qquad < \bigl( a+ \bigl( 4+ \Lambda T^{1/2} \bigr) t^{1/2} \bigr)
\biggl[ 1 + \int_{0}^t \frac{2 \Lambda}{1-\alpha}\exp
\biggl( \frac{2 \Lambda s}{1-\alpha} \biggr) \,ds \biggr]
+ \alpha g_{a}(t)
\\
&&\qquad = (1-\alpha) g_{a}(t) + \alpha g_{a}(t)
\\
&&\qquad =
g_{a}(t).
\end{eqnarray*}
The strict inequality remains true when $t=0$ since
${\mathbb E}[(X_{0})_{+}] \leq a < g_{a}(0)$.
Now, by the continuity of the paths of $Z^e$ and by\vspace*{1.5pt} the finiteness of\break
$\mathbb{E} [ \sup_{0 \leq t \leq T}
(Z_{t}^e)_{+}]$, we deduce that $\mathbb{E} [ \sup_{0 \leq s \leq t}
(Z_{s}^e)_{+}]$ is continuous in $t$. Therefore, if $R^{e} < T$, then
$\mathbb{E} [ \sup_{0 \leq s \leq R^e}
(Z_{s}^e)_{+}]$ must be equal to $g(R^e)$, but by the above
inequalities, this sounds as a contradiction.
By (\ref{integerrepresentation}),
this proves the announced bound.


\section{Existence and uniqueness in small time}\label{sectionExistenceanduniquenessinsmalltime}

The main result of this section is the following.

%
\begin{teo} \label{fixedpoint} Suppose
there exist $\beta,\varepsilon>0$ such that
$\mathbb{P}(X_{0}\in dx)\leq\beta(1-x)\,dx$ for any $x \in(1-\varepsilon,1]$ and that
the density of $X_{0}$ on the interval $(1-\varepsilon,1]$ is
differentiable at point 1. Then there exist constants $A_{1}\geq0$ and
$T_{1}\in(0,1]$, depending upon $\beta, \varepsilon, \alpha, \Lambda$
and $K$ only, such
that
$\Gamma(\mathcal{H}(T_{1},A_{1}) ) \subset\mathcal
{H}(T_{1},A_{1})$. Moreover, for all $e_{1},e_{2}\in\mathcal{H}(T_{1},A_{1})$,
\[
\bigl\llVert\Gamma(e_{1})-\Gamma(e_{2})\bigr\rrVert
_{\mathcal{H}(T_{1},A_{1})}\leq\tfrac{1}{2}\llVert e_{1}-e_{2}
\rrVert_{\mathcal{H}(T_{1},A_{1})}.
\]
Hence, there exists a unique fixed point of the restriction of $\Gamma
$ to $\mathcal{H}(T_{1},A_{1})$,
which provides a solution to (\ref{simplifiedeq}) according to
Definition~\ref{definitionsolution} up until time $T_{1}$ (such that
$[0, T_1] \ni t\mapsto\mathbb{E}(M_t)$ is in the space $\mathcal
{H}(T_{1},A_{1})$).
\end{teo}

\subsection{Representation of \texorpdfstring{$\Gamma$}{$Gamma$}}
Let $T> 0$. As a first step toward understanding the map $\Gamma$
defined above, we note that, given $e \in{\mathcal L}(T,A)$, using the
definitions we can write
\begin{eqnarray*}
\Gamma(e) (t) &=&\mathbb{E}\bigl(M_{t}^{e}\bigr)=
\mathbb{E} \biggl( \sum_{k \geq1} {\mathbh
1}_{[0,t]}\bigl(\tau_{k}^e\bigr) \biggr)
\\
&=& \sum_{k\geq1} \int_{0}^t
\mathbb{P} \bigl( \tau_{k+1}^e \in(s,t] |
\tau_{k}^e = s \bigr) \mathbb{P}\bigl(
\tau_{k}^e \in ds\bigr) +\mathbb{P}\bigl(
\tau_{1}^{e}\leq t\bigr),
\end{eqnarray*}
where $\mathbb{P}( \tau_{k}^e \in ds)$ is a convenient abuse of
notation for denoting the law of $\tau_{k}^e$ and
${\mathcal B}(\mathbb{R}) \ni \mathcal{A} \mapsto\mathbb{P}(\tau_{k+1}^e \in
\mathcal{A} |\tau_{k}^e = s)$ stands for the conditional law of
$\tau_{k+1}^e$ given $\tau_{k}^e=s$. Here, ${\mathcal B}(\mathbb
{R})$ is the Borel $\sigma$-algebra on $\mathbb{R}$. Moreover,
observing that the solution $X^e$ to
(\ref{gammadefinition2}) is a Markov process (which restarts from $0$
at time $\tau_{k}^e$ when $k \geq1$), we may write
%
%
\begin{eqnarray}\label{Markovproperty}
\Gamma(e) (t)&=&\mathbb{E}\bigl(M_{t}^{e}\bigr)
\nonumber\\[-8pt]\\[-8pt]\nonumber
&=& \sum_{k\geq1}\int_{0}^{t}
\mathbb{P} \bigl(\tau_{1}^{e^{\sharp_{s}}} \leq t-s |
X_{0}^{e^{\sharp_{s}}}=0 \bigr)\mathbb{P}\bigl(\tau_{k}^{e}
\in ds\bigr)+\mathbb{P}\bigl(\tau_{1}^{e} \leq t
\bigr),
\end{eqnarray}
where $e^{\sharp_{s}}$ stands for the mapping $([0,T-s] \ni t \mapsto
e(t+s)-e(s) ) \in{\mathcal L}(T-s,A)$.

With this decomposition it is clear that in order to analyze $\Gamma
(e)$, and more importantly the derivative of $\Gamma(e)$ [recall we
are looking for a fixed point in $\mathcal{H}(T, A)$], we must analyze
the densities of the first hitting times of a barrier by a \textit
{nonhomogeneous} diffusion process with a general Lipschitz drift
term. Indeed, formally taking the derivative with respect to $t$ in
(\ref{Markovproperty}) introduces terms involving the density of $\tau
_{1}^{e}$, where we recall that
\[
\tau_{1}^{e} = \inf\bigl\{t>0\dvtx  X_{t-}^e
\geq1\bigr\} = \inf\biggl\{t>0\dvtx  X_{0}+\int_{0}^{t}b
\bigl(X_{s}^{e}\bigr)\,ds+W_{t} \geq1 -\alpha e(t)
\biggr\}.
\]
The analysis of such densities is well known to be a difficult problem.
These problems remain even in the case where $b\equiv0$. However, the
fact that $e$ is continuously differentiable at least guarantees that
the densities exist. In the case $b\equiv0$, we refer to
[\citet{Peskir}, Theorem 14.4].
In the general case, existence of these densities will be guaranteed in
the next section by Lemma~\ref{lemkilledprocess1}.

\subsection{General bounds for the density of the first hitting time for a nonhomogeneous diffusion process}
\label{Generalbounds}

Fix $T>0$, and for $e\in\mathcal{C}^1[0,T]$ consider the stochastic
process $(\chi^e_{t})_{0 \leq t \leq T}$ which satisfies
%
%
\begin{equation}
\label{eq21102} d \chi^e_{t} = b\bigl(
\chi^e_{t}\bigr) \,dt + \alpha e'(t) \,dt +
dW_{t}, \qquad t \in[0,T], \chi^e_{0}
<1\mbox{ a.s.},
\end{equation}
together with the stopping time
\[
\tau^e:= \inf\bigl\{t \in[0,T]\dvtx  \chi^e_{t}
\geq1\bigr\}, \qquad(\inf\varnothing= \infty).
\]
Here, $\alpha\in(0,1)$ and the drift $b$ is globally Lipschitz,
exactly as above.

%
\begin{lemma}
\label{lemkilledprocess1}Let $e\in\mathcal{C}^1[0, T]$.
Suppose there exist $\beta,\varepsilon>0$ such that
$\mathbb{P}(\chi_{0}\in dx)\leq\beta(1-x)\,dx$ for any $x \in
(1-\varepsilon,1]$ and that
the density of $\chi_{0}$ on the interval $(1-\varepsilon,1]$ is
differentiable at point 1. Then:
\begin{longlist}[(iii)]
\item[(i)] For any $t \in(0,T]$, the law of
the diffusion $\chi^e_t$ killed at the threshold is absolutely
continuous with respect to the Lebesgue measure.

\item[(ii)] Denoting the density of $\chi^e_t$ killed at the threshold by
%
%
\begin{equation}
\label{eq21101} p_e(t,y):= \frac{1}{dy} \mathbb{P} \bigl(
\chi^e_{t} \in dy, t < \tau^e \bigr), \qquad t
\in[0,T], y \leq1,
\end{equation}
$p_e(t,y)$ is continuous in $(t, y)$ and continuously differentiable in
$y$ on
$(0,T] \times(-\infty,1]$ and admits Sobolev derivatives of order 1
in $t$ and of order 2 in $y$
in any $L^{\varsigma}$, $\varsigma\geq1$,
on any compact subset of $(0,T] \times(-\infty,1)$.
When $\chi_{0} \leq1-\varepsilon$ a.s.
it is actually continuous and continuously differentiable in $y$ on
any compact subset of $([0,T] \times(-\infty,1]) \setminus
(\{0\} \times(-\infty,1-\varepsilon])$.

\item[(iii)]
Almost everywhere on $(0,T] \times(-\infty,1)$,
$p_e$ satisfies
the Fokker--Planck equation:
%
%
\begin{equation}
\label{Fokker-Planck} \partial_{t} p_e(t,y) +
\partial_{y} \bigl[ \bigl( b(y) + \alpha e'(t) \bigr)
p_e(t,y) \bigr] - \tfrac{1}{2} \partial^2_{yy}
p_e(t,y) =0,
\end{equation}
with the Dirichlet boundary condition $p_e(t,1)=0$ and the
measure-valued initial condition
$p_e(0,y) \,dy = \mathbb{P}(\chi_{0} \in dy)$, $p_e(t,y)$ and $\partial
_{y}p_e(t,y)$ decaying to $0$ as
\mbox{$y \rightarrow- \infty$}.

\item[(iv)] The first hitting time, $\tau^e$ has a density on $[0,T]$, given by
%
%
\begin{equation}
\frac{d}{dt}\mathbb{P}\bigl(\tau^e \leq t\bigr)=-
\frac{1}{2}\partial_{y}p_e(t,1), \qquad t \in[0,T],
\label{derivofstoppingtime}
\end{equation}
the mapping $[0,T] \ni t \mapsto\partial_{y}p_{e}(t,1)$ being
continuous and its supremum norm being bounded in terms of $T$, $\alpha
$, $\|e'\|_{\infty, T}$, $\beta$ and $b$ only.
\end{longlist}
\end{lemma}

Lemma~\ref{lemkilledprocess1} is quite standard.
The analysis of the Green function of killed processes with smooth coefficients
may be found in [\citet{garronimenaldi}, Chapter VI].
The need for considering Sobolev derivatives follows from the fact that
$b$ is Lipschitz only. The argument to pass from the case $b$ smooth to
the case $b$ Lipschitz only is quite standard: it follows from
Calderon and Zygmund estimates, see [\citet{stroockvaradhan},
equation (0.4), Appendix~A],
that permit the control of the $L^{\varsigma}$ norm of the
second-order derivatives
on any compact subset of $(0,T] \times(-\infty,1)$.
A complete proof may be also found in the unpublished notes of
\citet{NOTES}.

When $\chi_0 = x_0$ for some deterministic $x_0<1$, the conditions of
the above lemma are certainly satisfied. Therefore, for $e\in\mathcal
{C}^1[0, T]$ it makes sense to consider the density $p_e(t, y)$, $t\in
(0, T]$, $y\leq1$ of the process killed at $1$ started at $x_0$. We
will write $p_e(t, y) = p^{x_0}_e(t, y)$ in this case.
The following two key results on $\partial_{y} p_{e}(t,1)$ are
standard adaptations of heat kernel estimates
[see, e.g., \citet{friedman}, Chapter~1] for killed processes.
The first one may be found in [\citet{garronimenaldi}, Chapter
VI, Theorem 1.10]
when $b$ is smooth and bounded. As explained in the beginning of
[\citet{garronimenaldi}, Chapter VI, Section~1.5] it remains true
when $b$ is Lipschitz continuous and bounded.
The argument for removing the boundedness assumption on $b$ is explained
in \citet{delaruemenozzi} in the case of a nonkilled process. As
shown in the unpublished
notes [\citet{NOTES}, Corollary 4.3], it can be adapted to the
current case.
The second result then follows from the so-called \textit{parametric}
perturbation argument following [\citet{friedman}, Chapter~1].
Again, the complete
proof can be found in the unpublished notes [\citet{NOTES},
Corollary 5.3].

\begin{prop}
\label{densityestimateprop}
Let $e\in\mathcal{C}^1[0, T]$. Then there exists a constant $\kappa
(T)$ (depending only on $T$ and the drift function $b$) which increases
with $T$ such that for all $x_0<1$,
\[
\bigl|\partial_yp_e^{x_0}(t, 1)\bigr| \leq\kappa(T)
\bigl(\bigl\|e'\bigr\|_{\infty, T} + 1\bigr)\frac{1}{t}\exp\biggl(-
\frac{(1-x_0)^2}{\kappa(T) t} \biggr)
\]
for all $t \leq\min\{[(\|e'\|_{\infty, T}+1)\kappa(T)]^{-2},T \}$.
In particular, $\kappa(T)$ is independent of~$e$.
\end{prop}

%
\begin{prop}
\label{densitydifferenceprop}
Let $e_1, e_2\in\mathcal{C}^1[0, T]$ and let $A = \max\{\|e'_1\|
_{\infty, T},\break \|e'_2\|_{\infty, T}\} $. Then there exists a constant
$\kappa(T)$ (depending only on $T$ and the drift function $b$) which
increases with $T$ such that for all $x_0<1$,
\begin{eqnarray*}
&& \bigl|\partial_yp_{e_{1}}^{x_{0}}(t, 1)-
\partial_yp_{e_{2}}^{x_{0}}(t, 1)\bigr|
\\
&&\qquad \leq
\kappa(T) (A + 1)\frac{1}{\sqrt{t}} \exp\biggl( -\frac{ (1-x_0
)^2}{\kappa(T)t} \biggr)
\bigl\|e'_1-e'_2\bigr\|_{\infty,t},
\end{eqnarray*}
for all $t \leq\min\{[(A+1)\kappa(T)]^{-2},T\}$. In particular,
$\kappa(T)$ is independent of $e_1$~and~$e_2$.
\end{prop}

\subsection{Application to \texorpdfstring{$\Gamma$}{$Gamma$}}
\label{subse53}
In this section, we return to the setting of Section~\ref
{Solutionasafixedpoint}, and apply the results of the previous
subsection to complete the proof of Theorem~\ref{fixedpoint}.

The first result ensures the differentiability of $\Gamma(e)$ whenever
$e\in\mathcal{L}(T, A)$, which is the first step in showing that
$\Gamma$ is stable on the space $\mathcal{H}(T, A)$ for some~$A$
(recall that $\mathcal{H}$ is simply a growth controlled subspace of
$\mathcal{L}$).

%
\begin{proposition}
\label{propdifferentiabilityGamma}
Let $e \in{\mathcal L}(T,A)$ and $X_{0}$ be such that there exist
$\beta,\varepsilon>0$ with
$\mathbb{P}(X_{0}\in dx)\leq\beta(1-x)\,dx$ for any $x \in(1-\varepsilon,1]$, and suppose that
the density of $X_{0}$ on the interval $(1-\varepsilon,1]$ is
differentiable at point 1.
Then the mapping $[0,T] \ni t \mapsto\Gamma(e)(t)$ is continuously
differentiable. Moreover,
%
%
\begin{eqnarray}\label{Markovproperty2}
\quad \frac{d}{dt} \bigl[ \Gamma(e) \bigr] (t) = - \int_{0}^{t}
\frac{1}{2}\partial_{y}p_{e}^{(0,s)}(t-s,1)
\frac{d}{ds} \bigl[ \Gamma(e) \bigr] (s) \,ds -\frac{1}{2}
\partial_{y}p_{e}(t,1),
\nonumber\\[-8pt]\\[-12pt]
\eqntext{t \in[0,T],}
\end{eqnarray}
where $p_{e}$ represents the density of the process $X^e$ killed
at $1$ and $p_{e}^{(0,s)}$ represents the density of the process
$X^{e^{\sharp_{s}}}$ killed
at $1$ with $X_{0}^{e^{\sharp_{s}}}=0$.
\end{proposition}

\begin{pf} We first check that $\Gamma(e)$ is Lipschitz continuous on $[0,T]$.
Considering a finite difference in (\ref{Markovproperty}) and using
(\ref{derivofstoppingtime}), we get, for $t,t+h \in[0,T]$,
%
%
\begin{eqnarray}\label{eq31010}
&& \Gamma(e) (t+h) - \Gamma(e) (t)\nonumber
\\
&&\qquad = \sum
_{k \geq1} \int_{t}^{t+h} \mathbb{P}
\bigl( \tau_{1}^{e^{\sharp_{s}}} \leq t+h-s | X^{e^{\sharp_{s}}}_{0}=
0 \bigr) \mathbb{P} \bigl(\tau_{k}^e \in ds\bigr)
\nonumber\\[-8pt]\\[-8pt]\nonumber
&&\quad\qquad{}{} - \frac{1}{2} \sum_{k \geq1}
\int_{0}^t \int_{t-s}^{t+h-s}
\partial_{y} p_{e}^{(0,s)}(r,1)\,dr \mathbb{P} \bigl(
\tau_{k}^e \in ds\bigr)
\\
&&\quad\qquad{} - \frac{1}{2} \int
_{t}^{t+h} \partial_{y} p_{e}(s,1) \,ds.\nonumber
\end{eqnarray}
By Lemma~\ref{lemkilledprocess1}(ii), we can handle the two last
terms in the above to find a constant $C>0$ (which depends on $e$)
such that
\begin{eqnarray*}
&& \Gamma(e) (t+h) - \Gamma(e) (t)
\\
&&\qquad \leq\sum_{k \geq1}
\int_{t}^{t+h} \mathbb{P} \bigl(
\tau_{1}^{e^{\sharp_{s}}} \leq t+h-s | X^{e^{\sharp_{s}}}_{0}=
0 \bigr) \mathbb{P} \bigl(\tau_{k}^e \in ds\bigr)
\\
&&\quad\qquad{}+ C h \bigl( 1 + \Gamma(e) (T) \bigr),
\end{eqnarray*}
the last term in the right-hand side being finite thanks to (\ref
{eq942}) and the argument following it.
Moreover, by (\ref{eq942}) and Gronwall's lemma, we deduce that
%
%
\begin{eqnarray}\label{eq31011}
&& \lim_{h \searrow0} \sup_{0 \leq s \leq T-h}
\mathbb{P} \bigl( \tau_{1}^{e^{\sharp_{s}}} \leq h |
X^{e^{\sharp_{s}}}_{0}=0 \bigr)
\nonumber
\\[-8pt]
\\[-8pt]
&&\qquad = \lim_{h \searrow0} \sup_{0 \leq s \leq T-h} \mathbb{P} \Bigl(
\sup_{0 \leq r \leq h} Z_{r}^{e^{\sharp_{s}}} \geq1 |
X^{e^{\sharp_{s}}}_{0}= 0 \Bigr) = 0,
\nonumber
\end{eqnarray}
where $Z^{e^{\sharp_{s}}}$ is
given by (\ref{integerrepresentation}).
Therefore, there exists a mapping $\eta\dvtx  \mathbb{R}_{+} \rightarrow
\mathbb{R}_{+}$ matching $0$ at $0$ and continuous at $0$ such that
\[
\Gamma(e) (t+h) - \Gamma(e) (t) \leq\eta(h) \bigl[ \Gamma(e)
(t+h) - \Gamma(e)
(t) \bigr] + C h \bigl( 1 + \Gamma(e) (T) \bigr).
\]
Choosing $h$ small enough, Lipschitz continuity easily follows.

As a consequence, we can divide both sides of (\ref{eq31010}) by $h$
and then let $h$ tend to~$0$. By (\ref{eq31011}), we have for a given
$t \in[0,T)$,
\begin{eqnarray*}
&&\lim_{h \searrow0} h^{-1} \sum
_{k \geq1} \int_{t}^{t+h} \mathbb{P}
\bigl( \tau_{1}^{e^{\sharp_{s}}} \leq t+h-s | X^{e^{\sharp_{s}}}_{0}=0
\bigr) \mathbb{P} \bigl(\tau_{k}^e \in ds\bigr)
\\
&&\qquad \leq\lim_{h \searrow0} \biggl[ \sup_{0 \leq s \leq T-h}
\mathbb{P} \bigl( \tau_{1}^{e^{\sharp_{s}}} \leq h |
X^{e^{\sharp_{s}}}_{0}=0 \bigr) \frac{\Gamma(e)(t+h)-\Gamma
(e)(t)}{h} \biggr] = 0.
\end{eqnarray*}
Handling the second term in (\ref{eq31010}) by Lemma
\ref{lemkilledprocess1} and using the Lebesgue dominated convergence theorem,
we deduce that
\begin{eqnarray}
\frac{d}{dt}\Gamma(e) (t) &=&-\sum_{k\geq1}\int
_{0}^{t}\frac{1}{2}\partial_{y}p_{e}^{(0,s)}(t-s,1)
\mathbb{P}\bigl(\tau^e_{k}\in ds\bigr)-\frac{1}{2}
\partial_{y}p_{e}(t,1).
\nonumber
\end{eqnarray}
By Lemma
\ref{lemkilledprocess1}, we know that $\partial
_{y}p_{e}^{(0,s)}(\cdot,1)$ and $\partial_{y} p_{e}(\cdot,1)$ are
continuous (in $t$). This proves that
$(d/dt)\Gamma(e)$ is continuous as well.

Formula (\ref{Markovproperty2}) then follows from the relationship
%
%
\begin{equation}
\label{eq5101} \Gamma(e) (t) = \sum_{k\geq1}\int
_{0}^t \mathbb{P}\bigl(\tau_{k}^e
\in ds\bigr), \qquad t \in[0,T].
\end{equation}\upqed
\end{pf}

The second idea is to show that the difference between the derivatives
of $\Gamma(e_1)$ and $\Gamma(e_2)$ is uniformly small in terms of the
distance between two functions $e_1$ and $e_2$ in the space $\mathcal
{H}(T, A)$ in small time.

%
\begin{prop} \label{differenceofderivatives} Let $T>0$ and $X_{0}$ be
such that there exist $\beta,\varepsilon>0$ with
$\mathbb{P}(X_{0}\in dx)\leq\beta(1-x)\,dx$ for any $x \in(1-\varepsilon,1]$, and suppose that
the density of $X_{0}$ on the interval $(1-\varepsilon,1]$ is
differentiable at point 1.

Suppose $e_{1},e_{2}\in\mathcal{H}(T, A)$ for some $A\geq0$. Then
there exists a constant $\kappa(T)$, independent
of $A$, $\beta$ and $\varepsilon$, and increasing in $T$,
and a constant $\tilde{\kappa}(T,\beta,\varepsilon)$, independent
of $A$ and increasing in $T$,
such that for any $e_{1},e_{2}\in\mathcal{H}(T, A)$,
\[
\sup_{0 \leq s\leq t}\biggl|\frac{d}{ds} \bigl[
\Gamma({e_{1}}) - \Gamma({e_{2}}) \bigr] (s)\biggr|
\leq(A+1) \tilde{\kappa}(T,\beta,\varepsilon) \sqrt{t} \bigl\|e'_{1}-e'_{2}
\bigr\|_{\infty,t},
\]
for $t \leq\min\{[(A+1)\kappa(T)]^{-2},T\}$.
\end{prop}

\begin{pf} We have by (\ref{Markovproperty2})
%
%
\begin{eqnarray}\label{eq4101}
&& \biggl|\frac{d}{dt} \bigl[
\Gamma(e_{1}) - \Gamma(e_{2}) \bigr](t) \biggr|\nonumber
\\
&&\qquad \leq
\frac{1}{2}\int_{-\infty}^{1}\bigl|\bigl[
\partial_{y}p_{e_{1}}^{x}-\partial_{y}p_{e_{2}}^{x}
\bigr](t,1)\bigr|\mathbb{P}(X_{0}\in dx)\nonumber
\\
&&\quad\qquad{}+\frac{1}{2} \int_{0}^{t}
\bigl|\bigl[ \partial_{y}p_{e_{1}}^{(0,s)} -
\partial_{y}p_{e_{2}}^{(0,s)} \bigr] (t-s,1)\bigr
|\frac{d}{ds} \Gamma(e_{1}) (s)
\\
&&\quad\qquad{}+\frac{1}{2}\int_{0}^{t}\bigl
|\partial_{y}p_{e_{2}}^{(0,s)}(t-s,1)\bigr
|\biggl|\frac{d}{ds} \bigl[ \Gamma(e_{1}) -
\Gamma(e_{2}) \bigr](s)\biggr| \,ds\nonumber
\\
&&\qquad:= \frac{1}{2} ( L_{1} + L_{2} + L_{3}).\nonumber
\end{eqnarray}

Suppose $t\leq T$ and $\sqrt{t}\leq[(A+1)\kappa(T)]^{-1}$, where
$\kappa(T)$ is as in Proposition~\ref{densitydifferenceprop}. The
value of $\kappa(T)$ will be allowed to increase when necessary below.
Considering
the first term only, we can use Proposition~\ref
{densitydifferenceprop} to see that
\begin{eqnarray*}
L_{1} &\leq&(A+1)\beta\kappa(T) \biggl(\int_{1- \varepsilon}^{1}
\frac{1}{\sqrt{t}} \exp\biggl( -\frac{(1-x)^{2}}{\kappa(T)t}
\biggr) (1-x)\,d{x} \biggr)
\bigl\|e'_{1}-e'_{2}
\bigr\|_{\infty,t}
\\
&&{} + (A+1) \kappa(T) \biggl(\int_{- \infty}^{1-\varepsilon}
\frac{1}{\sqrt{t}} \exp\biggl( -\frac{(1-x)^{2}}{\kappa(T)t}
\biggr) \mathbb{P}
(X_{0} \in dx) \biggr)\bigl\|e'_{1}-e'_{2}
\bigr\|_{\infty,t}.
\end{eqnarray*}

We deduce that there exists a constant $\tilde{\kappa}(T,\beta,\varepsilon)>0$, which is independent of $A$ and which is allowed to
increase as necessary from
line to line below, such that
%
%
\begin{eqnarray}
\label{eq4102} %
L_{1} &\leq&(A+1) \beta\kappa(T) \sqrt{t}
\biggl(\int_{0}^{\infty}z \exp\biggl( -
\frac{z^{2}}{\kappa(T)} \biggr) \,dz \biggr)\bigl\|e'_{1}-e'_{2}
\bigr\|_{\infty,t}
\nonumber
\\
&&{} + (A+1) \kappa(T) \frac{1}{\sqrt{t}} \exp\biggl( -\frac
{\varepsilon^{2}}{\kappa(T)t}
\biggr) \bigl\|e'_{1}-e'_{2}
\bigr\|_{\infty,t}
\\
&\leq&(A+1) \tilde{\kappa}(T,\beta,\varepsilon) \sqrt{t} \bigl\|e'_{1}-e'_{2}
\bigr\|_{\infty,t}.
\nonumber
\end{eqnarray}
We can then use Proposition~\ref{densitydifferenceprop}
again to see that
\[
L_{2} \leq(A+1) \kappa(T) \sup_{0 < s \leq t} \biggl[
s^{-1/2} \exp\biggl( - \frac{1}{\kappa(T) s} \biggr) \biggr]
\Gamma(e_{1}) (t) \bigl\|e'_{1}-e'_{2}
\bigr\|_{\infty,t}.
\]
By Proposition~\ref{stability} [since $e_1\in\mathcal{H}(T, A)$], we
deduce that
%
%
\begin{equation}
\label{eq4103} L_{2} \leq(A+1)\kappa(T) \sqrt{t}
\bigl\|e'_{1}-e'_{2}
\bigr\|_{\infty,t},
\end{equation}
where $\kappa(T)$ has been increased as necessary, and we have used
the elementary inequality $\exp(-1/v)\leq v$ for all $v\geq0$.
We finally turn to $L_{3}$ in (\ref{eq4101}). By Proposition~\ref
{densityestimateprop}, we have that
%
%
\begin{eqnarray}
\label{eq4105} %
\bigl|\partial_{y}p_{e_{2}}^{(0,s)}(t-s,1)
\bigr|&\leq&\kappa(T) (A+1)\frac{1}{(t-s)} \exp\biggl( -
\frac{1}{\kappa(T)(t-s)} \biggr)
\nonumber
\\[-8pt]
\\[-8pt]
\nonumber
& \leq&\kappa(T) (A+1),
\end{eqnarray}
again by increasing $\kappa(T)$.
Thus, from (\ref{eq4101}), (\ref{eq4102}), (\ref{eq4103}) and (\ref
{eq4105}), we deduce
\begin{eqnarray*}
\biggl|\frac{d}{dt} \bigl[\Gamma(e_{1}) -
\Gamma(e_{2}) \bigr](t) \biggr|&\leq&(A+1)\tilde{\kappa}(T,
\beta,\varepsilon) \sqrt{t} \bigl\|e'_{1}-e'_{2}
\bigr\|_{\infty,t}
\\
&&{}+ (A+1) \kappa(T)\int_{0}^{t}\biggl|
\frac{d}{ds} \bigl[\Gamma(e_{1}) -\Gamma(e_{2})
\bigr](s)\biggr| \,ds.
\end{eqnarray*}
By taking the supremum over all $s\leq t$ in the above, we have, for
$t\leq(2\kappa(T)(A+1))^{-1}$
[which actually follows from the aforementioned condition $t\leq
(\kappa(T)(A+1))^{-2}$ by assuming w.l.o.g. $\kappa(T) \geq2$],
\[
\sup_{0 \leq s\leq t}\biggl|\frac{d}{ds} \bigl[
\Gamma(e_{1})-\Gamma(e_{2}) \bigr](s)\biggr|\leq
2(A+1) \tilde{\kappa}(T,\beta,\varepsilon) \sqrt{t} \bigl\|e'_{1}-e'_{2}
\bigr\|_{\infty,t}.
\]\upqed
\end{pf}

We can then finally complete this section with the proof of Theorem
\ref{fixedpoint}.

\begin{pf*}{Proof of Theorem~\ref{fixedpoint}}
Choose
$A_{1}=2\sup_{0\leq t\leq1}|(d/dt)\Gamma(0)(t)|+1$.
Note that $A_{1}$ depends on $\beta$. Then choose $T_{1} \leq\min\{
[(A_1 +1)\kappa(1)]^{-2},1\}$
such that
%
%
\begin{equation}
\sqrt{T_{1}}\tilde{\kappa}(1,\beta,\varepsilon) (A_{1}+1)
\leq\tfrac{1}{4},\label{smallnessoftime}
\end{equation}
where $\kappa(1)$ and $\tilde{\kappa}(1,\beta,\varepsilon)$ are as in
Proposition~\ref{differenceofderivatives}.
By that result, if $e\in\mathcal{H}(T_{1},A_{1})$ then
\begin{eqnarray*}
&&\biggl|\frac{d}{dt}\Gamma(e) (t)\biggr|= \frac{d}{dt}
\Gamma(e) (t)\leq\sqrt{t}\tilde{\kappa}(T_{1},\beta,\varepsilon) (
A_{1} +1) A_{1}+\frac{d}{dt}\Gamma(0) (t)
\end{eqnarray*}
for all $t\leq\min\{[(A_{1}+1) \kappa(T_{1})]^{-2},T_{1}\}=T_{1}$.
By definition, we have $T_{1}\leq1$ so that $\kappa(T_{1}) \leq
\kappa(1)$ and $\tilde{\kappa}(T_{1},\beta,\varepsilon)\leq\tilde
{\kappa}(1,\beta,\varepsilon)$. Therefore,
\[
\frac{d}{dt}\Gamma(e) (t)\leq\sqrt{t}\tilde{\kappa}(1,\beta,\varepsilon)
(A_{1} +1) A_{1}+\frac{d}{dt}\Gamma(0) (t)
\]
for all $t \leq T_{1}$.
Hence, for all $t \leq T_{1}$
\begin{eqnarray*}
\frac{d}{dt}\Gamma(e) (t) &\leq&\frac{A_{1}}{2}+\sup
_{0\leq t\leq1} \biggl(\frac{d}{dt}\Gamma(0) (t) \biggr)\leq
A_{1}
\end{eqnarray*}
by (\ref{smallnessoftime}), so that $\Gamma(e)\in\mathcal{H}(T_{1},A_{1})$.

To prove that $\Gamma$ is a contraction on $\mathcal{H}(T_{1},A_{1})$,
first note that for $e\in\mathcal{H}(T_{1},A_{1})$
\[
\bigl\|e'\bigr\|_{\infty,T_{1}}\leq\|e\|_{\mathcal{H}(T_{1},A_{1})}\leq2
\bigl\|e'\bigr\|_{\infty,T_{1}}
\]
by the mean-value theorem, since $e(0)=0$ and $T_{1} \leq1$. Thus, for
any $e_{1},e_{2}\in\mathcal{H}(T_{1},A_{1})$
\begin{eqnarray*}
\bigl\llVert\Gamma(e_{1})-\Gamma(e_{2})\bigr\rrVert
_{\mathcal{H}(T_{1},A_{1})} &\leq&2\bigl\llVert\Gamma(e_{1})'-
\Gamma(e_{2})'\bigr\rrVert_{\infty,T_{1}}
\\
&\leq& 2\sqrt{T_{1}}\tilde{\kappa}(T_{1},\beta,\varepsilon)
(A_{1} +1) \bigl\|e'_{1}-e'_{2}
\bigr\|_{\infty,T_{1}}
\\
&\leq&\tfrac{1}{2}\|e_{1}-e_{2}
\|_{\mathcal{H}(T_{1},A_{1})},
\end{eqnarray*}
by our choice of $T_{1}$ and using Proposition~\ref{differenceofderivatives}
once more. Since $\mathcal{H}(T_{1},A_{1})$ is a closed subspace
of $\mathcal{C}^{1}[0,T]$ (a complete metric space), the existence
of a fixed point for $\Gamma$ follows from the Banach fixed-point
theorem.
\end{pf*}


\section{Long-time estimates}
\label{Long-timeestimates}
In order to extend the existence and uniqueness from small time to any
arbitrarily prescribed interval, we need an {a priori} bound
for the Lipschitz constant of $e\dvtx  t \mapsto\mathbb{E}(M_{t})$ on any
finite interval $[0,T]$, which is given by Theorem~\ref
{teogradientbd}. The purpose of this section is to prove this result.

As already mentioned, the key point is inequality (\ref{keyfact}).
Loosely, it says that,
in (\ref{Brunel}), the particles that are below
$1-dx$ at time $t$ receive a kick of order $\alpha\mathbb{P}(X_{t}
\in dx) < dx$.
In other words, only the particles close to $1$ can jump, which guarantees
some control on the continuity of $e$. Precisely,
Proposition~\ref{propholderbd1} gives a bound for the $1/2$-H\"older
constant of $e$.
Inequality (\ref{keyfact}) is proved by using {a priori}
heat kernel bounds when $\alpha$ is small enough, this restriction determining
the value of $\alpha_{0}$ in Theorem
\ref{teogradientbd}.
Once the $1/2$-H\"older constant of $e$ has been controlled, we provide in
Lemma~\ref{lemgradientbd41}
a H\"older estimate of the oscillation (in space) of
$p$ in the neighborhood of $1$. The proof is an adaptation of
\citet{krylovsafonov}. Finally,
in Proposition~\ref{propgradientbd5}, a barrier technique yields a
bound for the Lipschitz constant of $p$ in the neighborhood of $1$.

In the whole section, for a given initial condition
$X_{0}=x_{0} <1$, we thus assume that there exists a solution to (\ref
{simplifiedeq}) according to Definition~\ref{definitionsolution}, that
is, such that $e\dvtx  [0,T] \ni t
\mapsto\mathbb{E}(M_{t})$ is continuously differentiable.

\subsection{Reformulation of the equation and a priori bounds for the solution}

In the whole proof, we shall use a reformulated version
of (\ref{simplifiedeq}), in a similar way to Proposition~\ref
{stability} [see (\ref{integerrepresentation})]. Indeed,
given a solution $(X_{t},M_{t})_{0 \leq t \leq T}$ to (\ref
{simplifiedeq}) on some interval $[0,T]$ according to Definition~\ref
{definitionsolution}, we set
$Z_{t} = X_{t} + M_{t}$, $t \in[0,T]$. Then $(Z_{t})_{0 \leq t \leq T}$
has continuous paths and satisfies
%
%
\begin{equation}
\label{eqgradientbd2} Z_{t} = X_{0} + \int_{0}^t
b(X_{s}) \,ds + \alpha\mathbb{E}(M_{t}) + W_{t},
\qquad t \in[0,T],
\end{equation}
where
%
%
\begin{equation}
\label{eq791} M_{t} = \Bigl\lfloor\Bigl( \sup_{0 \leq s \leq t}
Z_{s} \Bigr)_{+} \Bigr\rfloor= \sup_{0 \leq s \leq t}
\bigl\lfloor(Z_{s})_{+} \bigr\rfloor.
\end{equation}

The following is easily proved by adapting the proof of Proposition
\ref{stability}:

%
\begin{lemma}
\label{lemgradientbd1}
There exists a constant $B(T,\alpha,b)$, only depending upon $T$,
$\alpha$, $b$ and nondecreasing in $\alpha$, such that
%
%
\begin{equation}
\label{eq491} \sup_{0 \leq t \leq T} e(t) = e(T) \leq\mathbb{E}
\Bigl[
\sup_{0 \leq t \leq T} (Z_{t})_{+} \Bigr] \leq B(T,
\alpha,b).
\end{equation}
A possible choice for $B$ is
\[
B(T,\alpha,b) = \frac{\mathbb{E}[(X_{0})_{+}] + 4 T^{1/2} + \Lambda
T}{1 - \alpha} \exp\biggl( \frac{2 \Lambda T}{1-\alpha} \biggr).
\]
\end{lemma}

\subsection{Local H\"older bound of the solution}
We now turn to the critical point of the proof. Indeed, in the next
subsection, we shall prove that, for $\alpha$ small enough, the function
$t \mapsto e(t) = \mathbb{E}(M_{t})$ generated by some solution to
(\ref{simplifiedeq}) according to Definition~\ref{definitionsolution}
(so that $e$ is continuously differentiable) satisfies an {a
priori} \mbox{$1/2$-}H\"older bound, with an explicit H\"older constant. This
acts as the keystone of the argument to extend the local existence and
uniqueness result into a global one. As a first step, the proof
consists of establishing a local H\"older bound for $e$ in the case
when the probability that the process $X$ lies in the neighborhood of
$1$ is not too large.

\begin{lemma}
\label{lemholderbd5}
Consider
a solution $(X_{t})_{0 \leq t \leq T}$ to (\ref{simplifiedeq}) on some
interval $[0,T]$, with $T>0$ and initial condition $X_{0}=x_{0}<1$.
Assume in addition that there exists some time $t_{0} \in[0,T]$ and
two constants $\varepsilon\in(0,1)$ and
$c \in(0,1/\alpha)$ such that
for any Borel subset $A \subset[1-\varepsilon,1]$,
%
%
\begin{equation}
\label{eqholderbd1} \mathbb{P} (X_{t_{0}} \in A ) \leq c | A
|,
\end{equation}
where $| A |$ stands for the Lebesgue measure of $A$. Then, with
\[
{\mathcal B}_0 = \frac{ \exp(2 \Lambda) [(8 + 5c + 8\varepsilon^{-1})
\Lambda
+ 4 (2+ c+ \varepsilon^{-1}) ]}{1- c \alpha},
\]
it holds that, for any
$h \in(0,1)$,
\[
\left.
\begin{array}{l} {\mathcal B}_{0} \exp(2\Lambda h)
h^{1/2} \leq\varepsilon/2
\\[2pt]
t_{0}+h \leq T
\end{array}
\right\} \Rightarrow e(t_{0}+h) -
e(t_{0}) \leq{\mathcal B}_{0} h^{1/2}.
\]
\end{lemma}

\begin{pf} By the Markov property, we can assume $t_{0}=0$, with $T$
being understood as $T-t_{0}$. Indeed, setting
%
%
\begin{equation}
\label{eq1292} X_{t}^{\sharp_{t_{0}}}:= X_{t_{0}+t}, \qquad t
\in[0,T-t_{0}],
\end{equation}
we observe that, for $t \in[0,T-t_{0}]$,
%
%
\begin{eqnarray}\label{eq1293}
X_{t}^{\sharp_{t_{0}}} &=& X_{t_{0}} + \int
_{0}^t b\bigl(X_{r}^{\sharp_{t_{0}}}
\bigr) \,dr + \alpha\mathbb{E}(M_{t+t_{0}}-M_{t_{0}})
\nonumber\\[-8pt]\\[-8pt]
&&{} +
W_{t+t_{0}} - W_{t_{0}} - (M_{t+t_{0}} - M_{t_{0}} ).\nonumber
\end{eqnarray}
Here $M_{t+t_{0}} - M_{t_{0}}$ represents the number of times the
process $X$ reaches $1$ within the interval \mbox{$(t_{0},t+t_{0}]$}.
Therefore, this also matches the number of times the process $X^{\sharp
_{t_{0}}}$ hits $1$ within the interval $(0,t]$, so that
${X}^{\sharp_{t_{0}}}$ indeed\vspace*{1pt} satisfies the nonlinear equation (\ref
{simplifiedeq}) on $[0,T-t_{0}]$, with ${X}^{\sharp
_{t_{0}}}_{0}=X_{t_{0}}$ as initial condition and with respect~to the
shifted Brownian motion
$({W}_{t}^{\sharp_{t_{0}}}:= W_{t_{0}+t} - W_{t_{0}})_{0 \leq t \leq
T-t_{0}}$. In what follows, $t_{0}$~is thus assumed to be zero, the new
$T$ standing for the previous $T-t_{0}$ and the new $X_{0}$ matching
the previous $X_{t_{0}}$ and thus satisfying (\ref{eqholderbd1}).

For a given $h\in(0, 1)$, such that $h \leq T$, and a given ${\mathcal
B}_{0} >0$ (the value of which will be fixed later), we then define the
deterministic hitting time:
\[
R = \inf\bigl\{t \in[0,h]\dvtx  \mathbb{E} ( M_{t}) = e(t) \geq{
\mathcal B}_{0} h^{1/2} \bigr\}.
\]
Following the proof of (\ref{eq942}) [see more specifically (\ref
{eq1291})], we have, for any $t \in[0,h \wedge R]$,
\begin{eqnarray*}
M_{t} &\leq& \sup_{0\leq s\leq t}(Z_{s})_{+}
\\
&\leq&(X_{0})_{+} + \Lambda\int_{0}^{t}
\bigl( 1 + (Z_{s})_{+} + M_{s} \bigr) \,ds +
\alpha e(t) + 2 \sup_{0 \leq s \leq t} | W_{s} |
\\
&\leq&(X_{0})_{+} + 2\Lambda\int_{0}^{t}
( 1 + M_{s} ) \,ds + \alpha{\mathcal B}_{0} h^{1/2}
+ 2 \sup_{0 \leq s \leq t} | W_{s} |
\\
&\leq&(X_{0})_{+} + 2 \Lambda h + 2\Lambda\int
_{0}^t M_{s} \,ds + \alpha{\mathcal
B}_{0} h^{1/2} + 2 \sup_{0 \leq s \leq t} |
W_{s} |,
\end{eqnarray*}
where we have used (\ref{eq791}) to pass from the first to the second
line. By Gronwall's lemma, we obtain
%
%
\begin{eqnarray}
\label{eq391} %
M_{t} &\leq&\exp(2\Lambda h) \Bigl[
(X_{0})_{+} + 2 \Lambda h + \alpha{\mathcal
B}_{0} h^{1/2} + 2 \sup_{0 \leq s \leq h} |
W_{s} |\Bigr]
\nonumber
\\[-8pt]
\\[-8pt]
\nonumber
&\leq&(X_{0})_{+} + \exp(2\Lambda h) \Bigl[ 4
\Lambda h + \alpha{\mathcal B}_{0} h^{1/2} + 2 \sup
_{0 \leq s \leq h} | W_{s} |\Bigr],
\end{eqnarray}
as $\exp(2 \Lambda h) \leq1+ 2 \Lambda h \exp(2 \Lambda h)$ and
$(X_{0})_{+} \leq1$.

Assume that ${\mathcal B}_{0} \exp(2\Lambda h) h^{1/2} \leq\varepsilon
/2 \leq1/2$. Then, by Doob's
$L^2$ inequality for martingales,
%
%
\begin{eqnarray}
\label{eq495} %
\sum_{k \geq2} \mathbb{P} (
M_{t} \geq k ) &\leq&\sum_{k \geq2}\mathbb{P}
\Bigl( \exp(2\Lambda h) \Bigl[ 4\Lambda h + 2 \sup_{0 \leq s \leq
h} |
W_{s} |\Bigr] \geq k - 3/2 \Bigr)
\nonumber
\\
&\leq&2 \exp(2\Lambda h ) \mathbb{E} \Bigl[ 4\Lambda h + 2 \sup
_{0 \leq s \leq h} | W_{s} |\Bigr]
\\
&\leq&\exp(2\Lambda h ) \bigl[ 8 \Lambda h + 8 h^{1/2} \bigr].
\nonumber
\end{eqnarray}
Moreover,
\begin{eqnarray*}
&&\mathbb{P}(M_{t} \geq1)
\\
&&\qquad \leq\mathbb{P} \Bigl( (X_{0})_{+} + \exp(2\Lambda h)
\Bigl[ 4 \Lambda h + \alpha{\mathcal B}_{0} h^{1/2} + 2 \sup
_{0 \leq s \leq h} | W_{s} |\Bigr] \geq1 \Bigr)
\\
&&\qquad \leq\mathbb{P} \Bigl( X_{0} \in[1-\varepsilon,1], X_{0} +
\exp(2\Lambda h) \Bigl[ 4 \Lambda h + \alpha{\mathcal B}_{0}
h^{1/2} + 2 \sup_{0 \leq s \leq h} | W_{s} |
\Bigr] \geq1 \Bigr)
\\
&&\quad\qquad{}+ \mathbb{P} \Bigl( \exp(2\Lambda h) \Bigl[ 4 \Lambda h + 2
\sup
_{0 \leq s \leq h} | W_{s} |\Bigr] \geq\varepsilon/2
\Bigr)
\\
&&\qquad:= I_{1} + I_{2},
\end{eqnarray*}
where we have used ${\mathcal B}_{0} \exp(2\Lambda h) h^{1/2} \leq
\varepsilon/2$ in the third line.

By Doob's $L^1$ maximal inequality, we deduce that
%
%
\begin{equation}
\label{eqholderbd4} I_{2} \leq2 \exp(2\Lambda h) \varepsilon^{-1}
\mathbb{E} \bigl[ 4 \Lambda h + 2 | W_{h} |\bigr] \leq\exp(2\Lambda h)
\varepsilon^{-1} \bigl[ 8 \Lambda h + 4 h^{1/2} \bigr].
\end{equation}

We now switch to $I_{1}$. By independence of $X_{0}$ and $(W_{s})_{0
\leq s \leq T}$ and by (\ref{eqholderbd1}),
\begin{eqnarray*}
I_{1} &\leq& c \int_{0}^{\varepsilon} \mathbb{P}
\Bigl( \exp(2\Lambda h) \Bigl[ 4 \Lambda h + \alpha{\mathcal B}_{0}
h^{1/2} + 2 \sup_{0 \leq s \leq h} | W_{s} |
\Bigr] \geq x \Bigr) \,dx
\\
&\leq& c \int_{0}^{+ \infty} \mathbb{P} \Bigl( \exp(2
\Lambda h) \Bigl[ 4 \Lambda h + \alpha{\mathcal B}_{0}
h^{1/2} + 2 \sup_{0 \leq s \leq h} | W_{s} |
\Bigr] \geq x \Bigr) \,dx
\\
&=& c \exp(2\Lambda h) \mathbb{E} \Bigl[ 4 \Lambda h + \alpha
{\mathcal
B}_{0} h^{1/2} + 2 \sup_{0 \leq s \leq h} |
W_{s} |\Bigr].
\end{eqnarray*}
By Doob's $L^2$ inequality,
\[
I_{1} \leq c \exp(2\Lambda h) \bigl[ 4 \Lambda h + \alpha{\mathcal
B}_{0} h^{1/2} + 4 h^{1/2} \bigr].
\]
Together with (\ref{eqholderbd4}), we deduce that
\[
\mathbb{P}(M_{t} \geq1) \leq\exp(2\Lambda h) \bigl[ 4 \bigl(c + 2
\varepsilon^{-1}\bigr) \Lambda h + 4 \bigl(c+\varepsilon^{-1}\bigr)
h^{1/2} + c \alpha{\mathcal B}_{0} h^{1/2} \bigr].
\]
From (\ref{eq495}), we finally obtain, for $t \leq R \wedge h$,
\begin{eqnarray*}
\mathbb{E} (M_{t}) &=& \sum_{k\geq1}
\mathbb{P}(M_t\geq k)
\\
&\leq&\exp(2\Lambda h) \bigl[ 4 \bigl(2 + c + 2 \varepsilon^{-1}\bigr)
\Lambda h + 4 \bigl( 2+ c+ \varepsilon^{-1}\bigr) h^{1/2} + c
\alpha{\mathcal B}_{0} h^{1/2} \bigr]
\\
&\leq&\exp(2\Lambda h) \bigl[ \bigl(8 + 5c + 8 \varepsilon^{-1}\bigr)
\Lambda h + 4 \bigl( 2+ c+\varepsilon^{-1}\bigr) h^{1/2} \bigr] +
c \alpha{\mathcal B}_{0} h^{1/2},
\end{eqnarray*}
provided ${\mathcal B}_{0} \exp(2\Lambda h) h^{1/2} \leq\varepsilon
/2\leq1/2$, which implies
\[
c \alpha{\mathcal B}_{0} \exp(2\Lambda h) h^{1/2} \leq c
\alpha{\mathcal B}_{0} h^{1/2} + c \Lambda h,
\]
using the fact that $\exp(2\Lambda h) \leq1 + 2\Lambda h \exp
(2\Lambda h)$. Therefore, if $R \leq h$, then we can choose
$t= R$ in the left-hand side above. By continuity of $e$ on $[0,T]$, it
then holds $e(R)={\mathcal B}_{0} h^{1/2}$, so that
\begin{eqnarray*}
(1- c \alpha) {\mathcal B}_{0} h^{1/2} &\leq&\exp(2 \Lambda
h) \bigl[ \bigl(8 + 5 c + 8\varepsilon^{-1}\bigr) \Lambda h + 4 \bigl
( 2+
c+\varepsilon^{-1}\bigr) h^{1/2} \bigr]
\\
&<& \exp(2 \Lambda) \bigl[ \bigl(8 + 5 c + 8\varepsilon^{-1}\bigr)
\Lambda
+ 4 \bigl( 2+ c+\varepsilon^{-1}\bigr) \bigr] h^{1/2},
\end{eqnarray*}
which is not possible when
\[
{\mathcal B}_{0} = \frac{\exp(2 \Lambda) [
(8 + 5 c + 8\varepsilon^{-1})
\Lambda
+ 4 ( 2+ c+\varepsilon^{-1}) ]
}{1- c \alpha}.
\]
Precisely, with ${\mathcal B}_{0}$ as above and ${\mathcal B}_{0}\exp
(2 \Lambda h) h^{1/2} \leq\varepsilon/2$ it cannot hold $R \leq h$.
\end{pf}

\subsection{Global H\"older bound}
In this subsection, we shall prove the following.

%
\begin{proposition}
\label{propholderbd1}
Let $\varepsilon\in(0,1)$. Then there exists a positive constant $\alpha
_{0} \in(0,1]$, only depending upon $\varepsilon$, $K$ and
$\Lambda$, such that: whenever $\alpha< \alpha_{0}$, there exists a
constant ${\mathcal B}$, only depending on $\alpha$, $\varepsilon$,
$K$ and $\Lambda$, such that, for all
positive times
$T>0$ and initial conditions $X_{0}=x_{0} \leq1 - \varepsilon$, any
solution to (\ref{simplifiedeq}) according to Definition~\ref
{definitionsolution}
satisfies
\[
\left.
\begin{array} {l} {\mathcal B} h^{1/2} \leq\varepsilon/2
\\[2pt]
t_{0}+h \leq T
\end{array}
\right\} \Rightarrow e(t_{0}+h) -
e(t_{0}) \leq{\mathcal B} h^{1/2},
\]
for any
$h \in(0,1)$ and $t_{0} \in[0,T]$.
Note that ${\mathcal B}$ above may differ from ${\mathcal B}_{0}$ in
the statement of Lemma~\ref{lemholderbd5}.
The constant $\alpha_{0}$ can be described as follows. Defining
$T_{0}$ as the largest time less than 1 such that
\[
(1 - \varepsilon) \exp(\Lambda T_{0}) \leq1 - 7\varepsilon/8, \qquad
\Lambda T_{0} \exp(\Lambda T_{0}) \leq\varepsilon/8,
\]
$\alpha_{0}$ can be chosen as the largest (positive) real satisfying
[with $B(T_0, \alpha_0, b)$ as in Lemma~\ref{lemgradientbd1}]
\begin{eqnarray*}
\alpha_{0} B(T_{0},\alpha_{0},b) &\leq&
\varepsilon/4,
\\
\alpha_{0} 2^{3/2} \bigl(c'
\bigr)^{3/2} \exp\bigl(-\tfrac12 \bigr) \bigl[ \varepsilon^{-1} +
B(T_{0},\alpha_{0},b) \bigr] &\leq&1,
\\
\alpha_{0} \bigl[ c' T_{0}^{-1/2} +
2^{3/2} \bigl(c'\bigr)^{3/2} \exp\bigl(-\tfrac12 \bigr)
B(T_{0},\alpha_{0},b) \bigr] &\leq&1.
\end{eqnarray*}
Here, the constant $c'$ is defined by the following property: $c'>0$,
depending on $K$ only, is such that for any diffusion process
$(U_{t})_{0 \leq t \leq1}$
satisfying
\[
dU_{t} = F(t,U_{t}) \,dt + dW_{t}, \qquad t \in
[0,1],
\]
where $U_{0}=0$ and $F\dvtx  [0,T] \times\mathbb{R} \rightarrow\mathbb
{R}$ is $K$-Lipschitz in $x$ such that $F(t,0)=0$ for any
$t \in[0,1]$, it holds that
\[
\frac{1}{dx}\mathbb{P} ( U_{t} \in dx ) \leq\frac{c'}{\sqrt{t}}
\exp\biggl( - \frac{x^2}{c' t } \biggr), \qquad x \in\mathbb{R},
t \in(0,1].
\]
\end{proposition}

The proof relies on the following.

\begin{lemma}
\label{lemholderbd2}
Given an initial condition $X_{0}=x_{0} \leq1-\varepsilon$, with
$\varepsilon\in(0,1)$,
and a solution $(X_{t})_{0 \leq t \leq T}$ to (\ref{simplifiedeq}) on
some interval $[0,T]$ according to Definition~\ref
{definitionsolution}, the random variable $X_{t}$ has a density on
$(-\infty,1]$, for any $t \in(0,T]$. Moreover, defining $T_{0}$ as in
the statement of Proposition~\ref{propholderbd1} and choosing $\alpha
\leq\alpha_{1}$ satisfying
\[
\alpha_{1} B(T_{0},\alpha_{1},b) \leq
\varepsilon/4,
\]
it holds, for $x \in[1-\varepsilon/4,1)$,
\begin{eqnarray*}
\frac{1}{dx}\mathbb{P} (X_{t} \in dx ) &\leq&2^{3/2}
\bigl(c'\bigr)^{3/2} \exp\biggl(-\frac12 \biggr) \bigl[
\varepsilon^{-1} + B(T_{0},\alpha,b) \bigr] \qquad\mbox{if } t
\leq T_{0},
\\
\frac{1}{dx} \mathbb{P} (X_{t} \in dx ) &\leq& c'
T_{0}^{-1/2} + 2^{3/2} \bigl(c'
\bigr)^{3/2} \exp\biggl(-\frac12 \biggr) B(T_{0},\alpha,b) \qquad
\mbox{if } t > T_{0},
\end{eqnarray*}
where the constant $c'$ is also as in the statement of Proposition~\ref
{propholderbd1}.
\end{lemma}

Before we prove Lemma~\ref{lemholderbd2}, we introduce some materials.
As usual, we set
$e(t) =\mathbb{E}(M_{t})$, for $t \in[0,T]$, the mapping $e$ being
assumed to be continuously differentiable on $[0,T]$.
Moreover, with $(X_{t})_{0\leq t \leq T}$, we associate the sequence of
hitting times
$(\tau_{k})_{k \geq0}$ given by (\ref{M}).
We then investigate the marginal distributions of $(X_{t})_{0\leq t
\leq T}$. Given a Borel subset
$A \subset(-\infty,1]$, we write in the same way as in the proof of
(\ref{Markovproperty})
%
%
\begin{eqnarray}
\label{eqholderbd11} %
\mathbb{P}(X_{t} \in A) &=& \mathbb{P} (
X_{t} \in A, \tau_{1} >t )
\nonumber
\\[-8pt]
\\[-8pt]
\nonumber
&&{}+ \sum_{k \geq1} \int_{0}^t
\mathbb{P} ( X_{t} \in A, \tau_{k+1} >t |
\tau_{k} = s ) \mathbb{P} (\tau_{k} \in ds),
\end{eqnarray}
where the notation $\mathbb{P} ( \cdot|\tau_{k} = s)$ stands
for the conditional law given $\tau_{k}=s$.
Following (\ref{eq1292}) and (\ref{eq1293}), we can shift the system
by length $s\in[0, T]$.
Precisely, we know that $(X_{r}^{\sharp_{s}}:=X_{s+r})_{0\leq r \leq
T-s}$ satisfies
%
%
\begin{equation}
\label{eqholderbd30} {X}_{r}^{\sharp_{s}} = X_{s} + \int
_{0}^r b \bigl( X_{u}^{\sharp_{s}}
\bigr) \,du + \alpha e^{\sharp_{s}}(r) + W_{s+r} - W_{s}
-M_{r}^{\sharp_{s}},
\end{equation}
with
\begin{eqnarray*}
e^{\sharp_{s}}(r)&:=& e(s+r) - e(s),\qquad
{M}_{r}^{\sharp_{s}}:= M_{s+r}- M_{s}
\quad\mbox{and}
\\
\tau_{k}^{\sharp_{s}}&:=& \inf\bigl\{ u > \tau_{k-1}^{\sharp_{s}}\dvtx
X_{s+u-} \geq1 \bigr\}
\end{eqnarray*}
for $k \geq1$,
($\tau_{0}^{\sharp_{s}}:= 0$).
Conditionally on $\tau_{k} = s$, the law of
$(X_{r}^{\sharp_{s}})_{0 \leq r \leq T-s}$ until $\tau_{1}^{\sharp
_{s}}$ coincides with the law of
$(\hat{Z}_{r}^{\sharp_{s},0})_{0 \leq r \leq T-s}$ until the first
time it reaches $1$,
where, for a given ${\mathcal F}_{0}$-measurable initial condition
$\zeta$ with values
in $(-\infty,1)$, $(\hat{Z}_{r}^{\sharp_{s},\zeta})_{0 \leq r \leq
T-s}$ stands for the
solution of the SDE:
%
%
\begin{equation}
\label{eqgradientbd100} \hat{Z}_{r}^{\sharp_{s},\zeta} = \zeta+
\int
_{0}^r b \bigl( \hat{Z}_{u}^{\sharp_{s},\zeta}
\bigr) \,du + \alpha e^{\sharp_{s}}(r) + W_{r}, \qquad r \in[0,T-s].
\end{equation}
Below, we will write $\hat{Z}^{\zeta}_r$ for $\hat{Z}_{r}^{\sharp
_{0},\zeta}$.
By (\ref{eqholderbd11}),
%
%
\begin{eqnarray}
\label{eqholderbd11b} %
\mathbb{P}(X_{t} \in A) &\leq&\mathbb{P}
\bigl( \hat{Z}_{t}^{X_{0}} \in A \bigr) + \sum
_{k \geq1} \int_{0}^t \mathbb{P}
\bigl(\hat{Z}_{t-s}^{\sharp_{s},0} \in A \bigr) \mathbb{P}(
\tau_{k} \in ds)
\nonumber
\\[-8pt]
\\[-8pt]
\nonumber
&=& \mathbb{P} \bigl( \hat{Z}_{t}^{X_{0}} \in A \bigr)
+ \int_{0}^t \mathbb{P} \bigl(
\hat{Z}_{t-s}^{\sharp_{s},0} \in A \bigr) e'(s) \,ds,
\end{eqnarray}
for any Borel set $A\subset(-\infty, 1]$, the passage from the first
to the second line following from (\ref{eq5101}).


\begin{pf*}{Proof of Lemma~\ref{lemholderbd2}}
Given an initial condition $x_{0} \in(-\infty,1-\varepsilon]$ for
$\varepsilon\in(0, 1)$,
we know from \citet{delaruemenozzi} that $\hat{Z}_{t}^{x_{0}}$
has a density for any $t \in(0,T]$
(and thus $\hat{Z}_{t-s}^{\sharp_{s},0}$ as well for $0 \leq s < t$).
From (\ref{eqholderbd11b}), we deduce that the law of $X_{t}$ has a
density on $(-\infty,1]$ since
$\mathbb{P}(X_{t} \in A)=0$ when $| A |= 0$, where $| A
|$ stands for the Lebesgue measure of $A$.
Moreover, there exists a constant $c' \geq1$, depending on $K$ only,
such that, for any $t \in[0,T \wedge1]$:
\begin{equation}
\label{eq496} \frac{1}{dx}\mathbb{P}\bigl(\hat{Z}_{t}^{x_{0}}
\in dx\bigr) \leq\frac{c'}{\sqrt{t}} \exp\biggl( - \frac{[x-
\vartheta_{t}^{x_{0}}]^2}{c' t} \biggr),
\end{equation}
where $\vartheta_{t}^{x_{0}}$ is the solution of the ODE:
%
%
\begin{equation}
\label{eqdefinitionvartheta} \frac{d}{dt} \vartheta_{t} = b (
\vartheta_{t} ) + \alpha e'(t), \qquad t \in[0,T],
\end{equation}
with $\vartheta_{0}^{x_{0}} = x_{0}$. Above, the function $[0,T] \ni
t \mapsto e(t)$ represents
$[0,T] \ni t \mapsto\mathbb{E}(M_t)$ given $X_0=x_0$, which means
that the initial condition $x_{0}$ of $X_{0}$ upon which $e$ depends is
fixed once and for all, independently of the initial condition of
$\vartheta$. In particular, as the initial condition of $\vartheta$
varies, the function $e$ does not.
We emphasize that $c'$ is independent of $e$ and can be taken to be
that defined in Proposition~\ref{propholderbd1}. Indeed, we can write
$\mathbb{P}(\hat{Z}_{t}^{x_{0}} \in dx)$ as $\mathbb{P}( \hat
{Z}_{t}^{x_{0}} - \vartheta_{t}^{x_{0}} \in d(x - \vartheta
_{t}^{x_{0}}))$, with
\begin{eqnarray*}
d \bigl( \hat{Z}_{t}^{x_{0}} - \vartheta_{t}^{x_{0}}
\bigr) &=& F \bigl(t, \hat{Z}_{t}^{x_{0}} -
\vartheta_{t}^{x_{0}} \bigr) \,dt + dW_{t}, \qquad t
\in[0,T], \hat{Z}_{0}^{x_{0}} - \vartheta_{0}^{x_{0}}
= 0;
\\
F(t,x) &=& b \bigl( x + \vartheta_{t}^{x_{0}} \bigr) - b \bigl(
\vartheta_{t}^{x_{0}} \bigr), \qquad t \in[0,T],  x \in
\mathbb{R}.
\end{eqnarray*}
We then notice that $F(t,\cdot)$ is $K$-Lipschitz continuous (since
$b$ is) and satisfies $F(t,0)=0$, so that, referring to \citet
{delaruemenozzi}, all the parameters involved in the definition of the
constant $c'$ are independent of $e$.
The fact that $c'$ is independent of $e$ is crucial. As a consequence,
we can bound $(1/dx)\mathbb{P}(\hat{Z}_{t-s}^{\sharp_{s},0} \in dx)$
in a similar way, that is, with the same constant $c'$ as in (\ref
{eq496}): for any
$0 \leq s < t \leq T$, with $t-s \leq1$,
%
%
\begin{equation}
\label{eq497} \frac{1}{dx} \mathbb{P}\bigl( \hat{Z}_{t-s}^{\sharp_{s},0}
\in dx\bigr) \leq\frac{c'}{\sqrt{t-s}} \exp\biggl( - \frac
{[x-\vartheta_{t-s}^{\sharp_{s},0}]^2}{c'(t-s)} \biggr),
\end{equation}
where $\vartheta^{\sharp_{s},0}$ is the solution of the ODE:
\[
\frac{d}{dt} \vartheta_{t}^{\sharp_{s}} = b \bigl(
\vartheta_{t}^{\sharp_{s}} \bigr) + \alpha\frac{d}{dt}e^{\sharp_{s}}(t),
\qquad t \in[0,T-s],
\]
with $\vartheta_{0}^{\sharp_{s}, 0} = 0$ as initial condition.

\textit{Bound of the density in small time.} Keep in mind that
$X_{0}=x_{0} \leq1-\varepsilon$. Therefore, by the comparison principle
for ODEs, $\vartheta^{x_{0}}_{t} \leq\vartheta^{1-\varepsilon}_{t}$
for any $t \in[0,T]$, so that by Gronwall's lemma
\[
\vartheta_{t}^{x_{0}} \leq\vartheta_{t}^{1-\varepsilon}
\leq\bigl( 1 - \varepsilon+ \Lambda T + \alpha e(T) \bigr) \exp
(\Lambda T).
\]
By Lemma~\ref{lemgradientbd1}, we know that $e(T) \leq B(T,\alpha,b)$, so that
%
%
\begin{equation}
\label{eq1295} \vartheta_{t}^{x_{0}} \leq\bigl( 1 - \varepsilon
+ \Lambda T + \alpha B(T,\alpha,b) \bigr) \exp(\Lambda T).
\end{equation}
Now choose $T_0$ as in Proposition~\ref{propholderbd1}, that is,
$T_{0} \leq1$ such that
\[
(1 - \varepsilon) \exp(\Lambda T_{0}) \leq1 - 7\varepsilon/8, \qquad
\Lambda T_{0} \exp(\Lambda T_{0}) \leq\varepsilon/8,
\]
and then take $\alpha_{1} \in(0,1)$ such that
\[
\alpha_{1} B(T_{0},\alpha_{1},b) \exp(\Lambda
T_{0}) \leq\varepsilon/4.
\]
Then, whenever $\alpha\leq\alpha_{1}$, it holds that
\[
\vartheta_{t}^{x_{0}} \leq1-\varepsilon/2, \qquad t \in
[0,T_{0} \wedge T].
\]
Therefore, for $x \geq1- \varepsilon/4$,
%
%
\begin{equation}
\label{eq5102} \exp\biggl( - \frac{[x-\vartheta_{t}^{x_{0}}]^2}{c'
t} \biggr) \leq\exp\biggl( -
\frac{\varepsilon^2}{16 c' t} \biggr), \qquad t \in[0,T_{0} \wedge T].
\end{equation}
Similarly,
\[
\vartheta_{t-s}^{\sharp_{s},0} \leq3\varepsilon/8 \leq3/8, \qquad
0\leq
s \leq t \leq T_{0} \wedge T.
\]
Indeed, $e^{\sharp_{s}}(T-s) \leq e(T)$ for $s \in[0,T]$, so that
(\ref{eq1295}) applies to
$\vartheta_{t-s}^{\sharp_{s},0}$
with $1-\varepsilon$ therein being replaced by $0$.
Therefore, for $x \geq1- \varepsilon/4$, it holds that $x-\vartheta
_{t-s}^{\sharp_{s},0} \geq3/4-3/8 = 3/8 \geq1/4$, so that
%
%
\begin{equation}
\label{eq5103} \qquad\exp\biggl( - \frac{[x-\vartheta_{t-s}^{\sharp
_{s},0}]^2}{c' (t-s)} \biggr) \leq\exp\biggl( -
\frac{1}{16 c' (t-s)} \biggr), \qquad0 \leq s < t \leq T_{0} \wedge T.
\end{equation}
In the end, for $x \in(1-\varepsilon/4,1)$ and $t \leq T_{0} \wedge T$,
we deduce
from (\ref{eqholderbd11b}), (\ref{eq496}), (\ref{eq497}), (\ref
{eq5102}), (\ref{eq5103})
and Lemma~\ref{lemgradientbd1} again, that
%
%
\begin{equation}
\label{eq1111} \qquad\frac{1}{dx} \mathbb{P}( X_{t} \in dx ) \leq
c' \varpi_{0} \bigl[ \varepsilon^{-1} + e(T
\wedge T_{0}) \bigr] \leq c' \varpi_{0} \bigl[
\varepsilon^{-1}+ B(T_{0},\alpha,b) \bigr],
\end{equation}
where
\begin{eqnarray*}
\varpi_{0} &=& \sup_{t >0} \biggl[ t^{-1/2}
\exp\biggl( - \frac{1}{16 c' t} \biggr) \biggr] = 4 \sqrt{c'}\sup
_{u >0} \bigl[ u \exp\bigl(-u^2\bigr) \bigr] =
2^{3/2} \sqrt{c'} \exp\biggl(-\frac{1}2 \biggr).
\end{eqnarray*}

\textit{Bound of the density in long time.} We now discuss what
happens for $T > T_{0}$ and
$t \in[T_{0},T]$.
Then
%
%
\begin{eqnarray}\label{eqpidefinition}
&& \frac{1}{dx} \mathbb{P}( X_{t} \in dx )\nonumber
\\
&&\qquad\leq\frac{1}{dx} \mathbb{P} \bigl(X_{t} \in dx,
\tau_{1}^{\sharp_{t-T_{0}}} \leq T_{0} \bigr) +
\frac{1}{dx} \mathbb{P} \bigl( X_{t} \in dx, \tau_{1}^{\sharp_{t-T_{0}}}
> T_{0} \bigr)
\\
&&\qquad = \pi_{1} + \pi_{2},\nonumber
\end{eqnarray}
with $\tau_{1}^{\sharp_{t-T_{0}}} = \inf\{u>0\dvtx  X_{t-T_{0}+u-} \geq
1\}
= \inf\{u>0\dvtx  X^{\sharp_{t-T_{0}}}_{u-} \geq1\}$. The above
expression says that we split the event ($X_{t}$ is in the neighborhood
of $x$) into two disjoint parts according to the fact that $X$ reaches
the threshold or not within the time
window $[t-T_{0},t]$.
We have chosen this interval to be of length
$T_0$ in order to apply the results in small time.

We first investigate $\pi_{2}$. The point is that, on the event that
$\tau_{1}^{\sharp_{t-T_{0}}} > T_{0}$ and within the time window
$[t-T_{0},t]$, $X$ behaves as a standard diffusion
process without any jumps, namely as a process with the same dynamics
as $\hat{Z}^{\sharp_{t-T_{0}},X_{t-T_{0}}}$.
Following~(\ref{eq496}), we then have
%
%
\begin{eqnarray}\label{eq499}
\pi_{2} &=& \frac{1}{dx} \mathbb{P} \bigl( \hat{Z}_{T_{0}}^{\sharp
_{t-T_{0}},X_{t-T_{0}}}
\in dx, \tau_{1}^{\sharp_{t-T_{0}}} > T_{0} \bigr)
\nonumber
\\
&\leq&\frac{1}{dx} \mathbb{P} \bigl( \hat{Z}_{T_{0}}^{\sharp
_{t-T_{0}},X_{t-T_{0}}}
\in dx \bigr)
\nonumber\\[-8pt]\\[-8pt]\nonumber
&\leq& \sup_{z \leq1} \frac{1}{dx} \mathbb{P}
\bigl(\hat{Z}_{T_{0}}^{\sharp_{t-T_{0}},z} \in dx \bigr)
\\
&\leq&  c' T_{0}^{-1/2}.\nonumber
\end{eqnarray}
We now turn to $\pi_{1}$.
Here, we write
\begin{eqnarray*}
\pi_1 &=& \frac{1}{dx} \mathbb{P} \bigl(X_{t} \in dx,
\tau_{1}^{\sharp_{t-T_{0}}} \leq T_{0} \bigr)
\\
&=&  \sum_{k\geq1} \frac{1}{dx}\mathbb{P} \bigl(X_{t} \in
dx, \tau_{k}^{\sharp_{t-T_{0}}} \leq T_{0} <
\tau_{k+1}^{\sharp_{t-T_{0}}} \bigr)
\\
&=& \sum_{k\geq1}\int_0^{T_0}
\frac{1}{dx} \mathbb{P} \bigl(X_{t} \in dx, T_{0} <
\tau_{k+1}^{\sharp_{t-T_{0}}}|\tau_{k}^{\sharp_{t-T_{0}}} =s
\bigr) \mathbb{P}\bigl(\tau_{k}^{\sharp_{t-T_{0}}} \in ds\bigr)
\\
&=& \sum_{k\geq1}\int_0^{T_0}
\frac{1}{dx}\mathbb{P} \bigl(\hat{Z}_{T_0 -s}^{\sharp
_{s+t-T_{0}},0}\in dx,
T_{0} < \tau_{k+1}^{\sharp_{t-T_{0}}} \bigr) \mathbb{P}\bigl(
\tau_{k}^{\sharp_{t-T_{0}}} \in ds\bigr),
\end{eqnarray*}
since on the event $\{\tau_{k}^{\sharp_{t-T_{0}}} \leq T_{0} < \tau
_{k+1}^{\sharp_{t-T_{0}}}\}$, given that the $k$th (and last) jump of
$X$ in the interval $[t-T_0, t]$ occurs at time $t - T_0 +s$ with $s\in
[0, T_0]$, we have that the process $X_r$ for $r\in[t - T_0 +s, t]$
coincides with the process $\hat{Z}_{u}^{\sharp_{s+t-T_{0}},0}$ for
$u\in[0, T_0 -s]$.
Thus,
%
%
\begin{eqnarray}
\label{eq499b} %
\pi_{1} &\leq&\sum
_{k \geq1}\int_{0}^{T_{0}}
\frac{1}{dx} \mathbb{P} \bigl( \hat{Z}_{T_0 -s}^{\sharp
_{s+t-T_{0}},0} \in dx
\bigr) \mathbb{P} \bigl( \tau_{k}^{\sharp_{t-T_{0}}} \in ds\bigr)
\nonumber
\\[-8pt]
\\[-8pt]
\nonumber
&=& \int_{0}^{T_{0}} \frac{1}{dx}
\mathbb{P} \bigl( \hat{Z}_{T_{0}-s}^{\sharp_{s+t-T_{0}},0} \in dx
\bigr)
e'(s+t-T_{0}) \,ds.
\end{eqnarray}
By (\ref{eq497}), we have
\begin{eqnarray*}
&&\int_0^{T_0}\frac{1}{dx} \mathbb{P} \bigl(
\hat{Z}_{T_{0}-s}^{\sharp_{s+t-T_{0}},0} \in dx \bigr)e'(s+t -
T_0)\,ds
\\
&&\qquad \leq\int_{0}^{T_{0}}
\frac{c'}{\sqrt{T_{0}-s}} \exp\biggl( - \frac{[x- {\vartheta
}^{\sharp_{s+t-T_{0}},0}_{T_{0}-s}]^2}{c' (T_{0}-s)} \biggr) e'(s+t-T_{0})
\,ds.
\end{eqnarray*}
Recalling that $e^{\sharp_{t-T_{0}}}(s)= \mathbb
{E}(M_{s+t-T_{0}}-M_{t-T_{0}})$, it is well seen that the mapping
$[0,T_{0}] \ni s \mapsto e^{\sharp_{t-T_{0}}}(s)$ satisfies
Lemma~\ref{lemgradientbd1}, that is,
\begin{eqnarray*}
\sup_{0 \leq s \leq T_{0}} e^{\sharp_{t-T_{0}}}(s) &=& \sup_{0 \leq s
\leq T_{0}}
\bigl[ e(s+t-T_{0}) - e(t-T_{0}) \bigr]
\\
&=& e(t) -
e(t-T_{0}) \leq B(T_{0},\alpha,b).
\end{eqnarray*}
Therefore, we can follow the same strategy as in short time; see (\ref{eq5103})
and (\ref{eq1111}). Indeed, for $\alpha\leq\alpha_{1}$, by the
choice of $T_0$ as before, it holds that
\[
\pi_{1} \leq c' \varpi_{0}
B(T_{0},\alpha,b),
\]
for $x\in[1-\varepsilon/4, 1)$.
Using (\ref{eq499}) and the above bound, we deduce that, for $t \in[T_{0},T]$,
\[
\frac{1}{dx} \mathbb{P} (X_{t} \in dx ) \leq c'
\bigl[ T_{0}^{-1/2} + \varpi_{0} B(T_{0},
\alpha,b) \bigr].
\]\upqed
\end{pf*}

\begin{pf*}{Proof of Proposition
\ref{propholderbd1}}
Proposition~\ref{propholderbd1} follows from the combination of
Lemmas~\ref{lemholderbd5} and~\ref{lemholderbd2}. Indeed, given $T_0$
and $\alpha_0$ as defined in Proposition~\ref{propholderbd1}, then by
Lemma~\ref{lemholderbd2} it follows that $\mathbb{P}(X_t\in A)
<(1/\alpha)|A|$ for any Borel subset $A\subset[1-\varepsilon/4, 1]$,
any $\alpha<\alpha_0$ and any $t\in[0, T]$. The result follows by
Lemma~\ref{lemholderbd5}, with $\mathcal{B}$ being given by $\mathcal
{B}_0\exp(2\Lambda)$ with $\varepsilon$ in $\mathcal{B}_0$ replaced
by $\varepsilon/4$.
\end{pf*}

\subsection{Estimate of the density of the killed process}
In light of the previous subsection,
for a solution $(X_{t})_{0 \leq t \leq T}$ to (\ref{simplifiedeq})
such that
the mapping $[0,T] \ni t \mapsto e(t)=\mathbb{E}(M_{t})$ is
continuously differentiable, we here investigate
\[
\frac{1}{dx} \mathbb{P} (X_{t} \in dx, t < \tau_{1}
), \qquad t \in[0,T], x \leq1,
\]
where $\tau= \inf\{t>0\dvtx  X_{t-}\geq1\}$ as usual. This is the density
of the killed process $(X_{t \wedge\tau_{1}})_{0\leq t \leq T}$,
which makes sense because of
Lemma~\ref{lemkilledprocess1}.

Here is the main result of this subsection.

\begin{lemma}
\label{lemgradientbd41}
Let $\varepsilon\in(0,1)$, $T >0$ and ${\mathcal B} >0$. Moreover, let
$(\chi_{t})_{0 \leq t \leq T}$ denote the solution to the SDE
\[
d \chi_{t} = b(\chi_{t}) \,dt + \alpha e'(t)
\,dt + dW_{t}, \qquad t \in[0,T]; \chi_{0}=x_0,
\]
for some continuously differentiable nondecreasing deterministic
mapping $[0,T] \ni t \mapsto e(t)$ satisfying
\[
e(0)=0, \qquad e(t) - e(s) \leq{\mathcal B} (t-s)^{1/2}, \qquad0 \leq
s\leq t\leq T.
\]
Then there exist two
positive constants $\mu_{T}$ and $\eta_{T}$, only depending upon $T$,
${\mathcal B}$, $\varepsilon$, $K$ and $\Lambda$, such that, for any
initial condition $x_{0} \leq1 - \varepsilon$,
%
%
\begin{equation}
\label{eq12927} p(t,y) \leq\mu_{T} (1-y)^{\eta_{T}}, \qquad t
\in[0,T], y \in[1-\varepsilon/4,1],
\end{equation}
where $p(t, y)$ denotes the density of $\chi_t$ killed at $1$ as in
(\ref{lemkilledprocess1}).
\end{lemma}

\begin{pf}
\textit{First step}.
The first step is to provide a probabilistic representation for $p$.
For a given $(T,x) \in(0,+\infty) \times(-\infty,1)$, we consider
the solution to the SDE:
%
%
\begin{equation}
\label{eq12925} dY_{t} = - \bigl[ b(Y_{t}) + \alpha
e'(T-t) \bigr] \,dt + dW_{t}, \qquad t \in[0,T],
Y_{0}=y,
\end{equation}
together with some stopping time $\rho\leq\rho_{0} \wedge T$, where
$\rho_{0}=\inf\{t \in[0,T]\dvtx  Y_{t} \geq1\}$
(with $\inf\varnothing= + \infty$). Then, by
Lemma~\ref{lemkilledprocess1} and
the It\^o--Krylov formula
[see \citet{krylov}, Chapter II, Section~10],
\begin{eqnarray*}
&&d \bigl( p (T-t,Y_{t}) \bigr)
\\
&&\qquad = - \partial_{t} p(T-t,Y_{t}) \,dt - \bigl[
b(Y_{t}) + \alpha e'(T-t) \bigr] \partial_{y}
p(T-t,Y_{t}) \,dt
\\
&&\quad\qquad{} + \tfrac{1}{2} \partial_{yy}^2
p(T-t,Y_{t}) \,dt
+ \partial_{y} p(T - t,Y_{t})
\,dW_{t}
\\
&&\qquad = b'(Y_{t}) p(T-t,Y_{t}) \,dt +
\partial_{y} p(T - t,Y_{t}) \,dW_{t},
\end{eqnarray*}
for $0 \leq t \leq\rho$.
Therefore, the Feynman--Kac formula yields
%
%
\begin{equation}
\label{eq22} p(T,y) = {\mathbb E} \biggl[ p(T-\rho,Y_{\rho})
{\mathbh
1}_{\{Y_{\rho} \neq1\}} \exp\biggl( - \int_{0}^{\rho}
b'(Y_{s}) \,ds \biggr) \Big| Y_{0} = y \biggr],
\end{equation}
the indicator function following from the Dirichlet boundary condition
satisfied by~$p(\cdot,1)$.

\textit{Second step}.
We now specify the choice of $\rho$.
Given some free parameters $L \geq1$ and $\delta\in(0,\varepsilon/4)$
such that $L \delta\leq\varepsilon/4$, we assume that the initial
condition $y$ in (\ref{eq12925}) is in $(1-\delta,1)$ and then
consider the stopping time
%
%
\begin{equation}
\label{eq12926a} \rho= \inf\bigl\{ t \in[0,T]\dvtx  Y_{t} \notin(1-L
\delta,1)\bigr\} \wedge\delta^2.
\end{equation}
Assume that $\delta^2 \leq T$. By (\ref{eq22}), we deduce that
%
%
\begin{equation}
\label{eq20} p(T,y) \leq\exp\bigl(K \delta^2\bigr) \bigl(1 - {
\mathbb P}(Y_{\rho}=1) \bigr) \sup_{(t,z) \in{\mathcal Q}(\delta,L)} p(t,z),
\end{equation}
with
\[
{\mathcal Q}(\delta,L) = \bigl\{(t,z) \in\bigl[T- \delta^2,T\bigr]
\times[1 - L \delta,1] \bigr\}.
\]
The point is then to give a lower bound for ${\mathbb P}(Y_{\rho}=1)$.
By assumption, we know that $e$ is $(1/2)$-H\"older continuous on
$[0,T]$. Therefore, since $Y_{0}=y \in(1-\delta,1)$, we have,
for any $t \in[0,\rho]$,
\[
Y_{t} \geq1-\delta- m \delta^2 - \alpha{\mathcal B}
\delta+ W_{t},
\]
with
%
%
\begin{equation}
\label{eq3173} m = \sup_{0 \leq z \leq1} \bigl| b(z) \bigr|.
\end{equation}
Therefore, for $m \delta\leq1$,
\[
Y_{t} \geq1- 2 \delta- \alpha{\mathcal B} \delta+
W_{t}, \qquad t \in[0,\rho],
\]
so that
%
%
\begin{equation}
\label{eq21} \qquad\{Y_{\rho}=1 \} \supset\Bigl\{ \sup_{0 \leq t \leq
\delta^2}
W_{t} > (2 +\alpha{\mathcal B}) \delta\Bigr\} \cap\Bigl\{ \inf
_{0 \leq t \leq\delta^2} W_{t} > (2 + \alpha{\mathcal B}- L) \delta
\Bigr\}.
\end{equation}
Choosing $L=3 + \alpha{\mathcal B}$ and applying a scaling argument,
we deduce that
%
%
\begin{eqnarray}
\label{eq18101} %
&& \mathbb{P} \Bigl( \Bigl\{ \sup_{0 \leq t \leq\delta^2}
W_{t} > (2+ \alpha{\mathcal B}) \delta\Bigr\} \cap\Bigl\{ \inf
_{0 \leq t \leq\delta^2} W_{t} > (2 + \alpha{\mathcal B} - L) \delta
\Bigr\} \Bigr)
\nonumber
\\
&&\qquad = \mathbb{P} \Bigl( \Bigl\{ \sup_{0 \leq t \leq1}
W_{t} > (2 + \alpha{\mathcal B}) \Bigr\} \cap\Bigl\{ \inf
_{0 \leq t \leq1} W_{t} > - 1 \Bigr\} \Bigr)
\\
&&\qquad =:c'' \in(0,1).\nonumber
\end{eqnarray}
We note that the above quantity $c''$ is independent of $\delta$ and
$T$. Moreover, we deduce from
(\ref{eq21}) that $\mathbb{P}(Y_{\rho}=1) \geq c''$ and, therefore,
from (\ref{eq20}) that
\[
p(T,y) \leq\bigl(1-c''\bigr) \exp\bigl(K
\delta^2\bigr) \sup_{z \in{\mathcal I}(L\delta)} \sup_{t \in[0,T]}
p(t,z),
\]
with ${\mathcal I}(r)=[1-r,1]$, for $r>0$. Choosing $\delta$ small
enough such that
$(1-c'') \exp(K \delta^2) \leq(1-c''/2)$, we obtain
\[
p(T,y) \leq\biggl(1- \frac{c''}{2} \biggr) \sup_{z \in{\mathcal
I}(L\delta)}
\sup_{t \in[0,T]} p(t,z), \qquad y \in{\mathcal I}(\delta).
\]
Modifying $c''$ if necessary ($c''$ being chosen as small as needed),
we can summarize the above inequality as follows: for $\delta\leq c''$,
%
%
\begin{equation}
\label{eq23} p(T,y) \leq\bigl(1-c''\bigr) \sup
_{ z \in{\mathcal I}(L \delta)} \sup_{t \in[0,T]} p(t,z), \qquad y
\in{
\mathcal I}(\delta).
\end{equation}
We now look at what happens when $T \leq\delta^2$ in (\ref{eq20}).
In this case, we can replace $\rho$ in the previous argument by $\rho
\wedge T$. Observing that $p(T-\rho\wedge T,Y_{\rho\wedge T})=0$ on
the event $\{\rho\geq T\} \cup\{Y_{\rho\wedge T}=1\}$ (since
$p(0,\cdot)=0$ on $[1-\varepsilon/4,1]$) and following~(\ref{eq20}), we
obtain, for $y \in{\mathcal I}(\delta)$,
%
%
\begin{equation}
\label{eq20prime} \qquad p(T,y) \leq\exp\bigl( K \delta^2\bigr) \bigl[1
- {
\mathbb P} \bigl(\{Y_{\rho\wedge T}=1\} \cup\{\rho\geq T\} \bigr)
\bigr] \sup
_{(t,z) \in{\mathcal Q}'(\delta,L)} p(t,z),
\end{equation}
with ${\mathcal Q}'(\delta,L) = \{(t,z) \in[0,T] \times[1 - L \delta,1] \}$. Now, the right-hand side of (\ref{eq21}) is included in the
event $\{Y_{\rho\wedge T}=1\} \cup\{\rho\geq T\}$ so that (\ref
{eq18101}) yields a lower bound for
${\mathbb P}(\{Y_{\rho\wedge T}=1\} \cup\{\rho\geq T\})$. Therefore,
we can repeat the previous arguments in order to prove that
(\ref{eq23}) also holds when $T \leq\delta^2$, which means that
(\ref{eq23}) holds true in both cases.

Therefore, by replacing $T$ by $t$ in the left-hand side in (\ref
{eq23}) and by letting $t$ vary within $[0,T]$, we have in any case,
\[
\sup_{y \in{\mathcal I}(\delta)}\sup_{t \in[0,T]} p(t,y) \leq
\bigl(1-c''\bigr) \sup_{z \in{\mathcal I}(L \delta)} \sup
_{t \in[0,T]} p(t,z).
\]
By induction, for any integer $n \geq1$ such that $L^n \delta\leq
r_{0}$, with $r_{0} = c'' \wedge(\varepsilon/4)$,
\[
\sup_{y \in{\mathcal I}(\delta)} \sup_{t \in[0,T]} p(t,y) \leq
\bigl(1-c''\bigr)^n \sup
_{z \in{\mathcal I}( L^n \delta)} \sup_{t \in[0,T]} p(t,z).
\]
Given $\delta\in(0, r_{0}/L)$, the maximal value for $n$ is
$n = \lfloor\ln[r_{0}/\delta]/\ln L \rfloor$. We deduce that, for
any $\delta\in(0, r_{0}/L)$,
%
%
\begin{equation}
\label{eq13915} \sup_{y \in{\mathcal I}(\delta)} \sup_{t \in[ 0,T]}
p(t,y) \leq\bigl(1-c''\bigr)^{( \ln[r_{0}/ \delta]/ \ln L) -1}
\sup
_{z \in{\mathcal I}(\varepsilon/4)}\sup_{t \in[0,T]} p(t,z).
\end{equation}
Following (\ref{eq496}), we know that
%
%
\begin{equation}
\label{eq18102} \sup_{z \in{\mathcal I}(\varepsilon/4)}\sup_{t \in
[0,T]} p(t,z)
\leq\sup_{z \in{\mathcal I}(\varepsilon/4)} \sup_{t \in[0,T]} \biggl[
\frac{c_{T}}{\sqrt{t}} \exp\biggl( - \frac{[z- \vartheta
_{t}^{x_{0}}]^2}{c_{T} t} \biggr) \biggr],
\end{equation}
for some constant $c_{T}$ only depending upon $T$ and $K$ and where
$(\vartheta_{t}^{x_{0}})_{0 \leq t \leq T}$ stands for the solution of
the ODE
\[
\frac{d \vartheta}{d t} = b(\vartheta_{t}) + \alpha e'(t),
\qquad t \in[0,T]; \vartheta_{0} = x_{0}.
\]
Pay attention that we here use the same notation as in
(\ref{eqdefinitionvartheta}) for the solution of the above ODE but
here $e(t)$ is not given as some $\mathbb{E}(M_{t})$. Actually, we
feel that there is no possible confusion here. Notice also that $e$ is
fixed and does not depend upon the initial condition $x_{0}$.

By the comparison principle for ODEs and then by Gronwall's lemma,
we deduce from the fact that $e$ is $(1/2)$-H\"older continuous that
\[
\vartheta_{t}^{x_{0}} \leq\vartheta_{t}^{1-\varepsilon}
\leq\bigl[ 1 - \varepsilon+ \Lambda t + {\mathcal B} t^{1/2} \bigr]
\exp(
\Lambda t), \qquad t \in[0,T].
\]
Using the above inequality, we can bound the right-hand side in (\ref
{eq18102}). Precisely, the above inequality says that the exponential
term in the supremum decays exponentially fast as $t$ tends to $0$ so
that the term inside the supremum can be bounded when $t$ is small;
when $t$ is bounded away from $0$, the term inside the supremum is
bounded by $c_{T}/\sqrt{t}$.
It is plain to deduce that
%
%
\begin{equation}
\label{eq12926b} \sup_{z \in{\mathcal I}(\varepsilon/4)}\sup_{t \in
[0,T]} p(t,z)
\leq c_{T},
\end{equation}
for a new value of $c_{T}$, possibly depending on $\varepsilon$ as well.
Therefore, for $\delta\in(0,r_{0}/L)$, (\ref{eq13915}) yields
\[
\sup_{y \in{\mathcal I}(\delta)} \sup_{t \in[ 0,T]} p(t,y) \leq
\frac{c_{T}}{(1-c'')} \biggl( \frac{\delta}{r_{0}} \biggr)^{\eta},
\]
with $\eta= - \ln(1-c'')/ \ln L $. This proves (\ref{eq12927}) for
$y \in(1-r_{0}/L,1)$. Note that $\eta$ is here independent of $T$,
contrary to what is indicated in the statement of Lemma~\ref
{lemgradientbd41}. However, we feel it is simpler to indicate $T$ in
$\eta_{T}$ as the constant ${\mathcal B}$ in the sequel will be chosen
in terms of $T$ thus making $\eta$ depend on $T$. Using (\ref
{eq12926b}), we can easily extend the bound to any $y \in(1-\varepsilon
/4,1)$ by modifying if necessary the parameters $\mu_{T}$ and $\eta
_{T}$ therein.
This completes the proof.
\end{pf}

\subsection{Bound for the gradient}
Here is the final step to complete the proof of Theorem~\ref{teogradientbd}.

%
\begin{proposition}
\label{propgradientbd5}
Let $\varepsilon\in(0,1)$, $T >0$ and ${\mathcal B} >0$. Moreover, let
$(\chi_{t})_{0 \leq t \leq T}$ denote the solution to the SDE
\[
d \chi_{t} = b(\chi_{t}) \,dt + \alpha e'(t)
\,dt + dW_{t}, \qquad t \in[0,T], \chi_{0}=x_0,
\]
for some continuously differentiable nondecreasing deterministic
mapping $[0,T] \ni t \mapsto e(t)$ satisfying
\[
e(0)=0;\qquad e(t) - e(s) \leq{\mathcal B} (t-s)^{1/2}, \qquad0 \leq
s\leq t\leq T.
\]
Then there exists a constant ${\mathcal M}_{T} >0$, only depending upon
$T$, ${\mathcal B}$, $\varepsilon$, $K$ and $\Lambda$, such that, for
any initial condition $x_{0} \leq1 - \varepsilon$ and any integer $n$
such that $n \geq\lceil4/\varepsilon\rceil$,
\[
\bigl|\partial_{y} p(t,1) \bigr|\leq\frac{{\mathcal M}_{T} n^{-\eta
_{T} }}{1 - \exp[-{\mathcal M}_{T}^{-1}
(1+\alpha C_{T}) n^{-1}]} (1+\alpha
C_{T}), \qquad t \in[0,T],
\]
where $p(t, y)$ is the density of $\chi_t$ killed at $1$ as in (\ref
{lemkilledprocess1}), $\eta_T$ is as in Lemma~\ref{lemgradientbd41}, and
\[
C_{T} = \sup_{0 \leq t \leq T} e'(t).
\]
\end{proposition}

\begin{pf}
We consider the barrier function
%
%
\begin{equation}
\label{eq51013} q(t,y) = \Theta\exp( Kt) \bigl[1 - \exp\bigl(
\gamma(y-1)
\bigr) \bigr], \qquad t \geq0, y \in\mathbb{R},
\end{equation}
where $\gamma$ and $\Theta$ are free nonnegative parameters. Then,
for $t>0$ and $y<1$,
\begin{eqnarray*}
&&\partial_{t} q(t,y) + \bigl( b(y) + \alpha e'(t)
\bigr) \partial_{y} q(t,y) - \tfrac{1}{2} \partial^2_{yy}
q(t,y)
\\
&&\qquad = \Theta\exp(Kt) \exp\bigl(\gamma(y-1)\bigr) \bigl( - \bigl
(b(y) + \alpha
e'(t)\bigr) \gamma+ \tfrac{1}{2} \gamma^2 \bigr)
+ K q(t,y).
\end{eqnarray*}
Keeping in mind that $\sup_{0 \leq t \leq T} e'(t) = C_{T}$
and choosing
%
%
\begin{equation}
\label{eq1710defgamma} \gamma= 2 \bigl( \max(m,1) + \alpha
C_{T}\bigr),
\end{equation}
where $m=\sup_{0\leq z\leq1}|b(z)|$ as before, we obtain, for $t
\in[0,T]$ and $y \in(0,1)$,
\[
- \bigl(b(y) + \alpha e'(t) \bigr) \gamma+ \tfrac{1}{2}
\gamma^2 \geq- 2 \bigl( \max(m,1) + \alpha C_{T}
\bigr)^2 + 2 \bigl( \max(m,1) + \alpha C_{T}
\bigr)^2 = 0.
\]
Thus, for $t \in[0,T]$ and $y \in(0,1)$,
\[
\partial_{t} q(t,y) + \bigl( b(y) + \alpha e'(t) \bigr)
\partial_{y} q(t,y) - \tfrac{1}{2} \partial^2_{yy}
q(t,y) \geq K q(t,y) \geq- b'(y) q(t,y),
\]
which reads
%
%
\begin{equation}
\label{eq26} \partial_{t} q(t,y) + \partial_{y} \bigl[
\bigl( b(y) + \alpha e'(t) \bigr) q(t,y) \bigr] - \tfrac{1}{2}
\partial^2_{yy} q(t,y) \geq0.
\end{equation}
For a given integer $n \geq\lceil4/\varepsilon\rceil$, we choose
$\Theta$ as the solution of
%
%
\begin{equation}
\label{eq51014} \Theta\biggl[ 1- \exp\biggl(- \frac{2 ( \max
(m,1) + \alpha C_{T})}{n} \biggr)
\biggr] = \mu_{T} n^{-\eta_{T}},
\end{equation}
with $\mu_{T}$ and $\eta_{T}$ as in the statement of
Lemma~\ref{lemgradientbd41}. Pay attention that the factor in the
left-hand side cannot be $0$ as $\max(m,1) >0$. Notice also
$q$ thus depends upon~$n$.
By Lemma~\ref{lemgradientbd41}, we deduce that
\[
q \biggl(t,1-\frac{1}{n} \biggr) \geq p \biggl(t,1-\frac{1}{n}
\biggr), \qquad0 \leq t \leq T.
\]
Now, we can apply the comparison principle for PDEs [see \citet
{lieberman}, Chapter IX, Theorem 9.7]. Indeed, we also observe that
$q(0,y) \geq p(0,y)=0$ for
$y \in[1-1/n,1]$ and
$q(t,1)=p(t,1) = 0$ for $t \in[0,T]$. Therefore, by (\ref{eq26}), we have
%
%
\begin{equation}
\label{eq51012} p(t,y) \leq q(t,y), \qquad t \in[0,T], y \in\biggl[1-
\frac{1}{n},1 \biggr].
\end{equation}
Since $p(t,1)=0=q(t, 1)$, we deduce
%
%
\begin{eqnarray}\label{eqgradientbd50}
\bigl|\partial_{y} p(t,1) \bigr| &\leq& \bigl|
\partial_{y} q(t,1) \bigr|
\nonumber\\[-8pt]\\[-8pt]\nonumber
&=& \frac{2 \mu_{T} ( \max(m,1) + \alpha
C_{T}) n^{-\eta_{T}}}{1- \exp[- 2 (\max(m,1) + \alpha C_{T})/n]}
\exp(Kt).
\end{eqnarray}\upqed
\end{pf}

We now complete the proof of Theorem~\ref{teogradientbd}.
We make use of Proposition~\ref{propdifferentiabilityGamma}. Recall
(\ref{Markovproperty2})
\[
e'(t) = - \int_{0}^{t}
\frac{1}{2}\partial_{y}p^{(0,s)}(t-s,1)
e'(s) \,ds -\frac{1}{2}\partial_{y}p(t,1), \qquad t
\in[0,T],
\]
where $p$ represents the density of the process $X$ killed
at $1$ and $p^{(0,s)}$ represents the density of the process
$X^{{\sharp_{s}}}$
driven by $e^{\sharp_{s}}=e(\cdot+s)-e(s)$ [see (\ref
{eqholderbd30})] killed
at $1$ with $X_{0}^{\sharp_{s}}=0$ as initial condition.

By Proposition~\ref{propholderbd1} and Lemma~\ref{lemgradientbd1}, we
know that, for a given $s \in[0,T)$ and for the prescribed values of
$\alpha$, the mapping
$[0,T-s] \ni r \mapsto e^{\sharp_{s}}(r)$ is $1/2$-H\"older
continuous, the H\"older constant only depending upon
$T$, $\alpha$, $\varepsilon$, $K$ and $\Lambda$ (Proposition~\ref
{propholderbd1} permits to bound the increments of $e^{\sharp_{s}}$ on
small intervals and Lemma~\ref{lemgradientbd1} gives a trivial bound
for the increments of $e^{\sharp_{s}}$ on large intervals). Therefore, by
Proposition~\ref{propgradientbd5}, we know that
%
%
\begin{eqnarray}\label{eq13912}
\bigl|\partial_{y}p^{(0,s)}(t-s,1) \bigr|\leq
\frac{{\mathcal M}_{T} n^{-\eta_{T} }}{1 - \exp[-{\mathcal
M}_{T}^{-1} (1+ \alpha C_{T}) n^{-1}]} (1+ \alpha C_{T}),
\nonumber\\[-12pt]\\[-8pt]
\eqntext{t \in[s,T],}
\end{eqnarray}
for $n \geq\lceil4/\varepsilon\rceil$ and for some constant ${\mathcal
M}_{T}$ only depending upon
$T$, $\alpha$, $\varepsilon$, $K$ and $\Lambda$. The same bound also
holds true for
$\partial_{y}p(t,1)$.

We deduce that, for any $t \in[0,T]$ and any $n$ such that $n \geq
\lceil4/\varepsilon\rceil$,
\[
e'(t) \leq\frac{ {\mathcal M}_{T} n^{-\eta_{T}}}{1- \exp[-
{\mathcal M}_{T}^{-1} (1+ \alpha C_{T}) n^{-1} ]}(1+ \alpha C_{T})
\frac{e(T)+1}{2}.
\]
By Lemma~\ref{lemgradientbd1}, we have a bound for $e(T)=\mathbb
{E}(M_{T})$, which means that we
can bound $(e(T)+1)/2$ in the right-hand side above by modifying the
constant ${\mathcal M}_{T}$.
Recalling
\[
C_{T} = \sup_{0 \leq t \leq T} e'(t),
\]
we deduce that
%
%
\begin{equation}
\label{eq51010} C_{T} \bigl( 1 - \exp\bigl[ - {\mathcal
M}_{T}^{-1} (1 + \alpha C_{T}) n^{-1}
\bigr] \bigr) \leq{\mathcal M}_{T} (1+ \alpha C_{T})
n^{-\eta_{T}}.
\end{equation}
Choosing $n$ large enough such that the right-hand side is less than
$(1+ \alpha C_{T})/2$ (so that $n$ depends on $T$) and multiplying by
$\alpha$, we get [since $\alpha\in(0,1)$]
\[
\frac{\alpha C_{T}}{2} \leq\frac{1}{2} + \alpha C_{T} \exp\bigl[-
{\mathcal M}_{T}^{-1} (1+\alpha C_{T})
n^{-1} \bigr].
\]
This shows that $\alpha C_{T}$ must be bounded in terms of ${\mathcal
M}_{T}$ and $n$. Precisely, we have
\[
\alpha C_{T} \leq1 + 2 \sup_{r \geq0} \bigl[ r \exp
\bigl[- {\mathcal M}_{T}^{-1} (1+r) n^{-1} \bigr]
\bigr]:=R < + \infty.
\]
By (\ref{eq51010}), we deduce that
\[
C_{T} \leq\sup_{0 \leq r \leq R} \biggl[\frac{ {\mathcal M}_{T} (1+
r) n^{-\eta_{T}}}{
1 - \exp[ - {\mathcal M}_{T}^{-1} (1 + r) n^{-1} ] }
\biggr],
\]
which is independent of $\alpha$ (for $\alpha\in(0,\alpha_{0}]$),
as required.
\qed


\section{Proof of Theorem \texorpdfstring{\protect\ref{solutionuptoT}}{2.4}}
\label{sectionExistenceanduniquenessforalltime}

In this section, we put everything together to arrive at our goal,
which is the proof of Theorem~\ref{solutionuptoT}. We first need the
following lemma, which is a corollary of Theorem~\ref{teogradientbd}.
The point is that the result will allow us to reapply the fixed-point
result on successive time intervals, since it guarantees that the
conditions of the fixed-point result are satisfied at the final point
of any interval on which we know there is a solution.

\begin{lem} \label{densitybound} For any $T>0$, initial condition
$X_{0}=x_{0} <1$, and $\alpha<\alpha_0$, where $\alpha_0=\alpha
_0(x_0)$ is as in Theorem~\ref{teogradientbd},
there exists a constant $C_{\mathrm{den}}(T)$ depending only on $T$, $x_0$,
$K$ and $\Lambda$
such that any solution to (\ref{simplifiedeq}) on $[0, T]$ satisfies
\[
\frac{1}{dy}\mathbb{P} (X_{t}\in dy )\leq C_{\mathrm{den}}(T)
(1-y),
\]
for all $y \in(1-\varepsilon/8,1)$ and $t \in[0,T]$, with $\varepsilon
=\min(1,1-x_{0})$.
\end{lem}

\begin{pf}
We assume that $(X_{t})_{0\leq t\leq T}$ is a
solution to (\ref{simplifiedeq}) with $X_{0}=x_0$ up until time $T$,
and set $e(t) = \mathbb{E}(M_t)$.
Following the notation of Section~\ref{sectionExistenceanduniquenessinsmalltime} (see also the last part of
the proof
of Theorem~\ref{teogradientbd}), for $y\leq1$
and $t\leq T$, let
\begin{eqnarray*}
p(t,y)&:=&\frac{1}{dy}\mathbb{P} (X_{t}\in dy,t <
\tau_{1} ),
\\
p^{(0,s)}(t,y)&:=& \frac{1}{dy}\mathbb{P} \bigl(X_{t}^{\sharp_{s}}
\in dy,t < \tau_{1}^{\sharp_{s}} | X_{0}^{\sharp_{s}}=0
\bigr).
\end{eqnarray*}
By Theorem~\ref{teogradientbd}, we know that $e$ is ${\mathcal
M}_{T}$-Lipschitz continuous, so that by (\ref{eq51012}),
\[
p(t,y) \leq q(t,y), \qquad t \in[0,T], y \in\biggl[1- \frac{1}{n},1
\biggr],
\]
where $n$ stands for $\lceil4/\varepsilon\rceil$ and $q$ is given
by (\ref{eq51013}), with $\gamma$ and $\Theta$ being fixed
by
(\ref{eq1710defgamma}) and (\ref{eq51014}),
with $C_{T} = {\mathcal M}_{T}$. By the specific form of $q$, this says
that there exists a constant
$C_{T}'$, depending only on $T$, $x_0$, $K$ and $\Lambda$, such that
\[
p(t,y) \leq C_{T}' (1-y), \qquad t \in[0,T], y \in
\biggl[1- \frac{\varepsilon}{8},1 \biggr],
\]
using the elementary inequality $1 - \exp(-x) \leq x$ for $x\in
\mathbb{R}$.
Clearly, the same argument applies to $p^{(0,s)}(t-s,y)$, that is,
\[
p^{(0,s)}(t-s,y) \leq C_{T}' (1-y), \qquad0
\leq s < t \leq T, y \in\biggl[1- \frac{\varepsilon}{8},1 \biggr].
\]
Now, following the proof of (\ref{eqholderbd11}), we get for $t \in
[0,T]$ and $y \in[1-\varepsilon/8,1]$,
%
%
\begin{eqnarray}\label{eq101112}
\frac{1}{dy} \mathbb{P}( X_{t} \in
dy) &=& p(t,y) + \int_{0}^t p^{(0,s)}(t-s,y)
e'(s) \,ds
\nonumber\\[-8pt]\\[-8pt]
&\leq& C_{T}' \bigl( 1+ e(T) \bigr)(1-y),\nonumber
\end{eqnarray}
where we use Lemma~\ref{lemkilledprocess1} for justifying the passage
to the density in (\ref{eqholderbd11}). By Lemma
\ref{lemgradientbd1}, this completes the proof.
\end{pf}

Finally, we can then prove the main result of the present paper:

\begin{pf*}{Proof of Theorem~\ref{solutionuptoT}}
We would like a solution up until fixed time $T>0$. The idea is to
iterate the fixed-point result (Theorem~\ref{fixedpoint}), which is
possible thanks to Lemma~\ref{densitybound}.
Indeed, by Theorem~\ref{fixedpoint}, we have
that there exists a solution to (\ref{simplifiedeq}) with $X_0 = x_0$
up until some small time $T_{1}>0$. By Lemma~\ref{densitybound},
we thus have that
%
%
\begin{equation}
\label{eqestiminduction} \frac{1}{dy}\mathbb{P}(X_{T_{1}}\in dy)\leq
C_{\mathrm{den}}(T_{1}) (1-y), \qquad y \in\biggl[1-
\frac{\varepsilon}{8},1 \biggr],
\end{equation}
where $\varepsilon= \min(1-x_{0},1)$. If $T_{1}\geq T$
we are done. If not, we have the above density bound
for $(1/dy)\mathbb{P}(X_{T_{1}}\in dy)$. We also know from (\ref
{eq101112}) and Lemma
\ref{lemkilledprocess1} that the density of $X_{T_{1}}$ is
differentiable at $y=1$.
Therefore, we can
apply Theorem~\ref{fixedpoint} again to see that there exists a
solution to (\ref{simplifiedeq}) on some interval $[T_{1},T_{1}+T_{2}]$
starting from $X_{T_{1}}$. As $T_{2}$ only depends upon $X_{T_{1}}$
through $\varepsilon$
(this is the statement of Theorem~\ref{fixedpoint})
and $C_{\mathrm{den}}(T_{1})$
and as these quantities can be bounded in terms of $T$, $\varepsilon$,
$K$, $\Lambda$ only, we then see
that
\[
T_{2} \geq\phi(T)
\]
for some constant $\phi(T)$ that refers to $T$, $\alpha$, $\varepsilon
$, $K$, $\Lambda$ only. Now we know that there exists a solution to
(\ref{simplifiedeq}) with $X_0 = x_{0}$ on $[0,T_{1}+T_{2}]$.
If $T_{1}+T_{2}>T$
we are done. If not, by Lemma~\ref{densitybound}
once again,
\[
\frac{1}{dy} \mathbb{P}(X_{T_{1}+T_{2}}\in dy)\leq C_{\mathrm{den}}(T_{1}+T_{2})
(1-y), \qquad y \in\biggl[1-\frac{\varepsilon}{8},1 \biggr],
\]
and we can then repeat the argument $n$ times to get a solution
up until time $T_{1}+\cdots+T_{n}$, where all $T_{k}\geq\phi(T)$
for $k\geq2$, that is, each time step is of size at least $\phi(T)$.
It is then clear that there exists $n\geq1$ such that $T_{1}+\cdots
+T_{n}\geq T$,
and so we are done for the existence of a solution.

Uniqueness of the solution proceeds in the same way. Given another
solution $(X_{t}',M_{t}')_{0\leq t \leq T}$ on the interval $[0,T]$ in
the sense of Definition~\ref{definitionsolution}, it must satisfy the
{a priori} estimates in the statements of
Theorem~\ref{teogradientbd} and
Lemmas~\ref{lemgradientbd1} and~\ref{densitybound}. In particular,
dividing the interval $[0,T]$ into subintervals of length $\phi(T)$
[except for the last interval the length of which might be less than
$\phi(T)$], with the same $\phi(T)$ as above, we can apply the
contraction property in Theorem~\ref{fixedpoint} on each subinterval
iteratively. Precisely, choosing $A_{1}$ accordingly in Theorem~\ref
{fixedpoint}, we prove by induction that the two solutions coincide on
$[0,\phi(T)]$, $[0,2\phi(T)]$, and so on.
\end{pf*}

\section*{Acknowledgement}
The authors would like to thank D. Talay for his involvement in a
number of fruitful discussions.




%

\printaddresses
\end{document}